\documentclass[a4paper,11pt]{article}

\usepackage[
    a4paper,
    top=1.05in,
    bottom=1.05in,
    left=1.00in,
    right=1.00in
]{geometry}
\usepackage[T1]{fontenc}

\usepackage{amsmath,amsthm,mathtools}
\usepackage{amssymb,mathrsfs}
\usepackage{newtxtext,newtxmath}
\usepackage{microtype}
\microtypesetup{protrusion=true,expansion=true}
\allowdisplaybreaks[4]

\usepackage{enumitem}
\setlist[enumerate]{leftmargin=2.5em,itemsep=0.25em,topsep=0.35em}
\setlist[itemize]{leftmargin=2.2em,itemsep=0.25em,topsep=0.35em}

\usepackage{titlesec}
\titleformat{\section}{\normalfont\Large\bfseries}{\thesection}{0.75em}{}
\titleformat{\subsection}{\normalfont\large\bfseries}{\thesubsection}{0.75em}{}
\titlespacing*{\section}{0pt}{2.6ex plus 0.8ex minus 0.2ex}{1.2ex plus 0.2ex}
\titlespacing*{\subsection}{0pt}{2.0ex plus 0.6ex minus 0.2ex}{0.9ex plus 0.2ex}

\usepackage{fancyhdr}
\fancypagestyle{plain}{%
  \fancyhf{}%
  \fancyfoot[C]{\thepage}%
}

\usepackage[dvipsnames]{xcolor}
\definecolor{linkblue}{RGB}{20,67,128}
\definecolor{citegreen}{RGB}{25,105,80}
\definecolor{urlpurple}{RGB}{100,55,125}
\usepackage[colorlinks=true,
            linkcolor=linkblue,
            citecolor=citegreen,
            urlcolor=urlpurple,
            plainpages=false]{hyperref}
\usepackage{bookmark}
\hypersetup{pdftitle={On the rigidity of generalized cylinders (V29)},
            pdfauthor={Hanbing Fang and Yu Li}}
\numberwithin{equation}{section}
\makeatletter
\renewcommand{\theequation}{\thesection.\ifnum\value{equation}<10 0\fi\arabic{equation}}
\makeatother

\newtheorem{theorem}{Theorem}[section]
\newtheorem{proposition}[theorem]{Proposition}
\newtheorem{lemma}[theorem]{Lemma}
\newtheorem{corollary}[theorem]{Corollary}
\theoremstyle{definition}
\newtheorem{definition}[theorem]{Definition}
\newtheorem{assumption}[theorem]{Assumption}
\theoremstyle{remark}
\newtheorem{remark}[theorem]{Remark}

\newcommand{\R}{\mathbb R}
\newcommand{\N}{\mathbb N}
\newcommand{\Rm}{\operatorname{Rm}}
\newcommand{\Ric}{\mathrm{Ric}}
\newcommand{\Tr}{\mathrm{Tr}}
\newcommand{\Id}{\operatorname{Id}}

\newcommand{\D}{\mathcal D}

\newcommand{\norm}[1]{\left\lVert #1\right\rVert}

\newcommand{\ISD}{\operatorname{ISD}}
\newcommand{\bnorm}{\norm{b}_{L^2(\mathrm{d}\mu_N)}}

\title{On the rigidity of generalized cylinders}
\author{Hanbing Fang\qquad and\qquad Yu Li}
\date{\today}

\newcommand{\anorm}{\norm{a}_{L^2(\mathrm d\bar\mu)}}

\begin{document}
\maketitle
\begin{abstract}
We prove pointed Gromov--Hausdorff rigidity of cylinders $N\times\mathbb R^d$
whenever $N$ is a compact Einstein manifold rigid as a Ricci shrinker.
\end{abstract}
\tableofcontents
\section{Introduction}
\label{sec:introduction}

Ricci shrinkers generate self-similar solutions to the Ricci flow and
play a central role in the analysis of finite-time
singularities~\cite{Perelman2002}.  Complete Ricci shrinkers in dimensions
at most three have been classified; see
Hamilton~\cite{Hamilton1995}, Perelman~\cite{Perelman2002},
Naber~\cite{Naber2010},
Ni--Wallach~\cite{NiWallach2008}, and
Cao--Chen--Zhu~\cite{CaoChenZhu2008}.  In higher dimensions, a complete
classification may be out of reach.

Instead of seeking a complete classification, another approach is to
study normalized Ricci shrinkers of a fixed dimension as a moduli space
equipped with the pointed
Gromov--Hausdorff topology, choosing minimum points of their potentials
as base points.  Weak compactness under an entropy lower bound provides
a framework for this study~\cite{LiLiWang2021}.  From this viewpoint,
the rigidity problem asks whether sufficient pointed closeness to a
model forces an isometry of the underlying normalized shrinkers.

For compact shrinkers, rigidity is closely related to the local
deformation theory of the soliton equation.  Podest\`a--Spiro developed
a finite-dimensional analytic description of the local moduli
space~\cite{PodestaSpiro2013}.  Kr\"oncke used obstructions to
infinitesimal solitonic deformations to prove rigidity of complex
projective spaces of even complex dimension~\cite{Kroencke2015}.
Li--Zhang subsequently proved Gromov--Hausdorff rigidity of complex
projective spaces in every complex dimension~\cite{LiZhangCP2022}.
Sun--Zhu proved rigidity of $\mathbb S^2\times\mathbb S^2$ and of
$\mathbb S^2\times N$ under spectral assumptions on the compact
Einstein factor $N$, by identifying a third-order
obstruction~\cite{SunZhu2025}.  Further results on obstructed Einstein
and solitonic deformations of compact symmetric spaces were obtained
by Batat--Hall--Murphy--Waldron~\cite{BatatHallMurphyWaldron2024} and
Hall--Schwahn--Semmelmann~\cite{HallSchwahnSemmelmann2025}.

Noncompact shrinkers are essential in the singularity analysis of Ricci
flow: even flows
on closed manifolds can develop singularities modeled on noncompact
spaces.  In this setting, a natural rigidity question asks whether
closeness to a model on a sufficiently large bounded region forces a
global isometry.

For noncompact shrinkers, local deformation theory alone does not yield
a general rigidity criterion. The difficulty is that pointed closeness
leaves the geometry outside a compact set uncontrolled. In particular,
it does not provide a global diffeomorphism identifying a nearby
shrinker with the reference manifold, as emphasized in
\cite{ColdingMinicozzi2025,LiZhang2023}. This motivates studying models
with well-understood geometry at infinity, for which comparison maps
can be constructed on successively larger regions.

A particularly important class consists of generalized cylinders.
We call the product of a compact Ricci shrinker and the Gaussian
shrinker $(\mathbb R^d,g_E,|\vec x|^2/4)$ a generalized cylinder.
Li--Wang proved pointed Gromov--Hausdorff rigidity of the round cylinder
\mbox{$\mathbb S^{n-1}\times\R$}~\cite{LiWangCylinder}.
Colding--Minicozzi established rigidity for all round cylinders and
developed the gauge, almost-splitting, contraction, and extension
arguments used to compare shrinkers on expanding
regions~\cite{ColdingMinicozzi2025}.  Li--Zhang extended this framework
to cylinders with compact $H$-stable Einstein factors, allowing nontrivial infinitesimal
deformations of the compact factor under an obstruction of order
at most three~\cite{LiZhang2023}.
More recently, Bamler--Lai~\cite{BamlerLai2026} gave another proof of the
rigidity of round cylinders among Ricci shrinkers using their PDE--ODI principle.

In this paper, we prove that a cylinder with compact Einstein
factor is rigid whenever its compact factor is rigid as a Ricci
shrinker. This condition is necessary: nearby nonisometric compact
Ricci shrinkers give rise to nearby nonisometric cylinders.

Let $(N^m,g_N)$ be smooth, closed, connected, and positive-dimensional,
with $\Ric_{g_N}=\frac{g_N}{2}$. Define
\begin{equation}
 \Theta_N=\log\operatorname{Vol}_{g_N}(N)-\frac m2\log(4\pi e).
 \label{eq:intro-normalization}
\end{equation}
For $d\ge1$ set $n=m+d$ and
\[
 (\bar M,\bar g)=(N\times\R^d,g_N+g_E),\qquad
 \bar f=\Theta_N+\frac m2+\frac{|\vec{x}|^2}{4}.
\]
\begin{theorem}[Pointed rigidity]
\label{thm:intro-rigidity}
Assume that $(N,g_N)$ is rigid among normalized compact Ricci
shrinkers. Fix $p_N\in N$ and
$\bar p=(p_N,0)$. There exists $\epsilon>0$ such that any connected
complete normalized Ricci shrinker $(M^n,g,f,p)$ with
$p\in\operatorname*{argmin}_Mf$ and
\[
 d_{\mathrm{pGH}}((M,g,p),(\bar M,\bar g,\bar p))<\epsilon
\]
is isometric to $(\bar M,\bar g)$.
\end{theorem}

The same conclusion holds for a fixed finite quotient by a free
isometric action; see Theorem~\ref{thm:quotient-rigidity}.
Examples of rigid compact factors include the round spheres
and spherical space forms, complex projective spaces,
and the rigid spherical products established in
\cite{Kroencke2015,LiZhangCP2022,SunZhu2025}.
For the linear stability of compact symmetric spaces, see Cao--He~\cite{CaoHe2015}.

The central analytic ingredient is a rigidity inequality: the size of
a small deformation is bounded by a power of its error in satisfying
the shrinker equation. The error is measured by
$\Phi(g,f)=\frac12g-\Ric_g-\nabla_g^2f$, together with the error in
the chosen gauge. Schematically, the inequality takes the form
\[
 (\text{size of the deformation})^\rho
 \le C\,(\text{shrinker defect}).
\]
Here size and error are measured in Gaussian weighted norms after
fixing coordinates, normalization, and center. In particular, a
deformation satisfying the shrinker equation and the gauge condition
must vanish.

This estimate improves comparisons on bounded regions and thereby
drives the global proof. Pointed closeness first supplies a comparison
on a sufficiently large region. Localization and gauge fixing produce
an approximate shrinker whose Gaussian weighted errors are small.
The rigidity inequality improves the interior closeness, and the
contraction--extension argument enlarges the comparison region.
Iteration gives comparisons on exhausting regions and ultimately a
global isometry.

Linear estimates alone do not give this improvement because the
linearized equation has a kernel. Its directions consist of
deformations inherited from $N$ and deformations varying along the
Euclidean factor; the orthogonal complement is controlled by the
linearized operator. Colding--Minicozzi~\cite{ColdingMinicozzi2025}
established round-cylinder rigidity using quadratic rigidity in a
suitable gauge and the geometric machinery for extending comparisons.
Li--Zhang~\cite{LiZhang2023} combined obstructions on the compact factor of order at
most three with quadratic control of kernel directions depending on the
Euclidean variables under their spectral stability hypotheses. Here we
control every kernel direction depending nontrivially on the Euclidean
variables, including additional resonant modes, and allow arbitrary finite
obstruction order on $N$. No tensor nonresonance or
linear-stability assumption is imposed on $N$. The new analytic input is
the rigidity inequality for the full kernel and an estimate controlling
the interaction of these two kinds of obstruction. The geometric
contraction--extension framework is inherited from those works.

Write the product kernel as $K=K_N\oplus K_+$, where $K_N$ consists
of the directions inherited from $N$ and $K_+$ of the positive Hermite-degree
components, as described in Section~\ref{sec:local-analysis}. Let
$b\in K_N$ and $a\in K_+$. Write $|b|$ for the $L^2$ norm with respect
to the normalized volume measure
$\operatorname{Vol}_{g_N}(N)^{-1}\mathrm dV_{g_N}$.
The family $X_b$ is obtained in Subsection~\ref{sec:local-kuranishi}
by taking the product of a family on $N$ with the Gaussian shrinker;
its analytic reduced entropy is $e(b)=\mathcal W(X_b)$.
Rigidity makes $0$ an isolated
zero of $\mathrm De$, so the Lojasiewicz distance inequality gives
\begin{equation}
 |\mathrm De(b)|\ge c_K|b|^\ell
 \quad\text{for all sufficiently small }b\in K_N,
 \qquad \ell<\infty.
 \label{eq:intro-obstruction}
\end{equation}
Here $\ell$ is a finite integer, with $\ell=1$ when $K_N=0$.
Proposition~\ref{prop:compact-isolation} relates this estimate to
compact rigidity, using the characterization in~\cite{LiZhang2023}.

The positive-degree directions admit a single quadratic test. Write
$\bar X=(\bar g,\bar f)$ and let $\mathrm d\bar\mu$ be its normalized
Gaussian weighted measure. The metric--density chart $\mathcal C$ of
Subsection~\ref{sec:kerspace} records $X=(g,f)=\mathcal C(H)$, with
$H=(h,\sigma)$, through
$g=\bar g+h$ and $\mathrm d\mu_X=(1+\sigma)\,\mathrm d\bar\mu$.
Set $\iota(h)=(h,0)$ and $v_*=(|\vec{x}|^2-2d)g_N$. If $h_j$ is
the Hermite-degree-$j$ component of $h\in K$, then
Theorem~\ref{thm:radial-obstruction} gives
\begin{equation}
 \left\langle\mathrm D^2(\Phi\circ\mathcal C)_0
 [\iota(h),\iota(h)],v_*\right\rangle_{L^2}
 =-\sum_{j\ge1}j^2\|h_j\|_{L^2}^2,
 \qquad h=\sum_jh_j\in K.
 \label{eq:intro-radial}
\end{equation}
All positive-degree contributions have the same sign, so this radial
test detects every direction in $K_+$, including
additional resonances. The degree-zero component is instead controlled
by the obstruction estimate on $N$ in~\eqref{eq:intro-obstruction}.

To combine these estimates, we move the radial test along the family
$X_b$. Applying the reduced Hessian to this direction gives
an error bounded by $|\mathrm De(b)|$. Radial testing and
differentiation in the $K_N$ directions then control $|a|^2$ and
$|\mathrm De(b)|$ together, with a cubic error in $a$ that can be
absorbed independently of $\ell$. Consequently, a normal reduction of
fixed order suffices: we retain the family $X_b$ exactly
in $b$, solve the normal equation through quadratic order in $a$,
and retain the reduced entropy through cubic order in $a$. The normal
residual is cubic in $a$, so quadratic control of $a$ suffices for
every finite obstruction order on $N$.

For the precise conclusion, suppose $H=(h,\sigma)$ is sufficiently
small in the bounded $C^2$ norm and its density is normalized and
centered, meaning
$\int_{\bar M}\sigma\,\mathrm d\bar\mu=0$ and
$\int_{\bar M}\vec{x}\,\sigma\,\mathrm d\bar\mu=0$.
Let $b$ and $a$ be the $L^2$ projections of $h$ onto $K_N$ and $K_+$,
respectively, and put $\Xi=H-\hat H_{b,a}$, where $\hat H_{b,a}$ is the cutoff
family constructed in Section~\ref{sec:quantitative}. All Sobolev
norms below and the $L^2$ norm of $a$ use $\mathrm d\bar\mu$.
For $X=\mathcal C(H)$, set
$\D(X)=\|\Phi(X)\|_{L^2}+\|\delta_{\bar X}h\|_{H^1}$.
Theorem~\ref{thm:local-rigidity-input} proves
\begin{equation}
 \anorm^2+|b|^\ell+\|\Xi\|_{H^{2}}
 \le C\D(X),\qquad
 \|H\|_{H^{2}}\le C\D(X)^{1/\rho},
 \quad \rho=\max\{2,\ell\}.
 \label{eq:intro-local}
\end{equation}

On the weighted divergence-free slice, the divergence term vanishes.
Theorem~\ref{thm:local-rigidity-input} also supplies the higher-regularity
estimates used to convert this weighted control into improved
geometric comparisons. Section~\ref{sec:global-rigidity} completes the
contraction--extension argument and carries it out equivariantly for
finite quotients.

The paper is organized as follows. Section~\ref{sec:local-analysis}
sets up the weighted gauge and product kernel.
Section~\ref{sec:quantitative} constructs the family $X_b$ and the
normal reduction through second order and establishes the remainder
estimates. Section~\ref{sec:Einstein-rigidity} proves the quadratic
obstruction on $K_+$ and the rigidity inequality for cylinders with
compact Einstein factors.
Section~\ref{sec:global-rigidity} proves global rigidity and treats finite
quotients.

\paragraph{Acknowledgements.}
Hanbing Fang would like to express gratitude to the University of Science and Technology of China (USTC) for its warm hospitality and generous support. This work was completed during Hanbing Fang’s visit to USTC. Hanbing Fang is supported by Hong Kong RGC grants No. 14300623
and No. 14304225. Yu Li is supported by the
National Key Research and Development Program of China (Grant No.~2025YFA1018200),
NSFC-12522105, and research funds from the University of Science and Technology of
China and the Chinese Academy of Sciences.

\paragraph{AI declaration.}
ChatGPT (ChatGPT-6 Astra) was used to help identify the quadratic
obstruction identity in Theorem~\ref{thm:radial-obstruction}, expand
selected technical arguments, and improve the presentation.
The proof framework and initial proof sketches were developed by the
authors, who take full responsibility for the mathematical content
and correctness.

\section{Weighted analysis on generalized cylinders}
\label{sec:local-analysis}
\subsection{Weighted operators and variation formulas}\label{sec:WM}
Let $(M^n,g)$ be a complete Riemannian manifold and let
$f\in C^\infty(M)$. We write $X=(g,f)$, with \emph{tensor part} $X_{\mathrm t}=g$
and \emph{scalar part} $X_{\mathrm s}=f$. Likewise, for any tensor--scalar pair
$H=(h,k)$, we write $H_{\mathrm t}=h$ and $H_{\mathrm s}=k$. Set
\[
 \mathrm d\mu_X:=(4\pi)^{-\frac{n}{2}}e^{-f}\;\mathrm dV_g,
 \qquad \mathcal V(X):=\int_M\mathrm d\mu_X,
\]
where $\mathrm dV_g$ is the Riemannian volume form. The pair $X$ is
\emph{normalized} if $\mathcal V(X)=1$.
\phantomsection\label{def:weighted-operators}
We write $\Delta_f=\Delta-\nabla_{\nabla f}$ for the weighted
Laplacian, acting on functions or tensors, and
$\mathcal L_f=\Delta_f+2\Rm_g$ for the tensor stability operator,
where $(\Rm_g h)_{ij}=R_{ikjl}h^{kl}$.
The divergence and weighted divergence are
\[
 (\delta_gh)_j=\nabla^ih_{ij},\qquad
 \delta_X:=\delta_{g,f},\qquad
 (\delta_Xh)_j=(\delta_gh)_j-h_{ij}\nabla^if.
\]
The formal adjoint of $\delta_X$ with respect to $\mathrm d\mu_X$ is
\[
 \delta_X^*\omega=\delta_g^*\omega=-\frac{1}{2}L_{\omega^\sharp}g.
\]
Here $L_V$ is the Lie derivative in the direction $V$, and
$\omega^\sharp$ is the vector field dual to $\omega$ with respect to $g$.

\begin{definition}
\label{def:weighted-quantities}
For $X=(g,f)$ and a variation $H=(H_{\mathrm t},H_{\mathrm s})$ of $X$, define
\[
 \Phi(X):=\frac{1}{2}g-\Ric_g-\nabla_g^2f,
 \qquad q_X(H):=\frac{1}{2}\Tr_gH_{\mathrm t}-H_{\mathrm s}.
\]
\end{definition}

If $\Phi(X)=0$, then $(M,g,f)$ is a \emph{Ricci shrinker}.
The weighted entropy and its scalar Euler--Lagrange expression are
\begin{align*}
 \mathcal W(X)&:=\int_M(\mathrm R_g+|\nabla f|_g^2+f-n)\;\mathrm d\mu_X,\\
 w(X)&:=2\Delta_gf-|\nabla f|_g^2+\mathrm R_g+f-n.
\end{align*}
We write $\boldsymbol{\mu}(g,1)$ for Perelman's entropy at scale one;
see \cite{Perelman2002} and \cite[Section~5]{LiWangHeatKernel}.

The following identities follow from the standard metric and entropy
variation formulas; see \cite{Perelman2002} and
\cite[Section~2]{CaoZhu2012}.
\begin{lemma}
\label{lem:pair-first-variation}
For a pair $X=(g,f)$ and a variation $J=(h,k)$, one has
\begin{align*}
 \mathrm D\Phi_X[J]
 &=\frac{1}{2}\mathcal L_fh+\delta_X^*\delta_Xh+\nabla_g^2q_X(J)
       +\frac{1}{2}\bigl(\Phi(X)\circ h+h\circ\Phi(X)\bigr),\\
 \mathrm Dw_X[J]
 &=-(2\Delta_f+1)q_X(J)+\delta_X\delta_Xh
       +\langle\Phi(X),h\rangle_g.
\end{align*}
Here compositions of symmetric tensors use $g$ to raise the contracted
index.
\end{lemma}

On the reference product, unqualified derivatives, tensor norms, and
integrals use $\bar g$ and $\mathrm d\bar\mu$. In particular,
\[
 \|T\|_{H^{j}}^2:=\sum_{i=0}^j\int_{\bar M}
             |\bar\nabla^iT|^2\;\mathrm d\bar\mu,
 \qquad [T]_j:=\sum_{i=0}^j|\bar\nabla^iT|.
\]
That is, $H^j$ denotes the weighted Sobolev space.
For a pair the squared norm is the
sum of the two squared component norms, and $[H]_j=[H_{\mathrm t}]_j+[H_{\mathrm s}]_j$.
The corresponding norms on $N$ use $g_N$ and $\mathrm d\mu_N$.
The corresponding bounded $C^{j,\theta}$ norms use the fixed product
charts, with one fixed $\theta\in(0,1)$. Write $|b|=\bnorm$ and
$|a|=\anorm$ for the finite-dimensional parameter norms. Their dual
norms are used for derivatives of reduced functions.

\subsection{Ricci shrinkers and generalized cylinders}\label{sec:RSI}

For a connected compact normalized Ricci shrinker $(N^m,g_N,f_N)$, set
$X_N=(g_N,f_N)$. We have
\begin{equation}
 \Ric_{g_N}+\nabla_{g_N}^2f_N=\frac{1}{2}g_N,
 \qquad
 (4\pi)^{-\frac{m}{2}}\int_N e^{-f_N}\,\mathrm{d}V_{g_N}=1,
 \qquad
 \Theta_N:=\boldsymbol{\mu}(g_N,1).
 \label{eq:compact-shrinker}
\end{equation}
Under this normalization, the potential function $f_N$ satisfies
\begin{equation}
 \mathrm{R}_{g_N}+|\nabla f_N|^2-f_N=-\Theta_N.
 \label{eq:compact-potential-identity}
\end{equation}
In particular, if $\Ric_{g_N}=\frac{g_N}{2}$, then $f_N$ is constant and
$\Theta_N=f_N-m/2$.
For $d=n-m\ge1$, the corresponding generalized cylinder is
\begin{equation*}
 \bar M^n=N^m\times\R^d,
 \qquad
 \bar g=g_N+g_E,
 \qquad
 \bar f=f_N+f_E,
 \qquad
 f_E=\frac{|\vec{x}|^2}{4}.
\end{equation*}
We write $\bar X:=(\bar g,\bar f)$.
Here $\vec{x}=(x_1,\ldots,x_d)$ denotes the standard coordinate on $\R^d$,
and we usually write $p$ for a point of $N$.  Set
$\langle \vec{x}\rangle:=1+|\vec{x}|$.  We use the probability
measures
\begin{equation*}
 \begin{gathered}
  \mathrm{d}\mu_N:=\mathrm d\mu_{X_N}=(4\pi)^{-\frac{m}{2}}e^{-f_N}\,\mathrm{d}V_{g_N},
  \qquad \mathrm{d}\gamma_d:=(4\pi)^{-\frac{d}{2}}e^{-\frac{|\vec{x}|^2}{4}}\,\mathrm{d}\vec{x},\\
  \mathrm{d}\bar\mu:=\mathrm d\mu_{\bar X}=\mathrm{d}\mu_N\mathrm{d}\gamma_d
  =(4\pi)^{-\frac{n}{2}}e^{-\bar f}\,\mathrm{d}V_{\bar g}.
 \end{gathered}
\end{equation*}

Tracing \eqref{eq:compact-shrinker}, using \eqref{eq:compact-potential-identity}, and integrating with respect to
$\mathrm{d}\mu_N$ gives
\begin{equation}
 \mathrm{R}_{g_N}+\Delta f_N=\frac{m}{2},
 \qquad
 \Delta_{f_N}f_N=\frac{m}{2}+\Theta_N-f_N,
 \qquad
 \int_Nf_N\,\mathrm{d}\mu_N=\frac{m}{2}+\Theta_N.
 \label{eq:weighted-basic-identities}
\end{equation}
We also use the shrinker identities:
\begin{equation}
 \delta_{X_N}\Ric_{g_N}=0,
 \qquad
 (\Delta_{f_N}+2\Rm_{g_N})\Ric_{g_N}=\Ric_{g_N},
 \label{eq:shrinker-Ricci-identities}
\end{equation}
\begin{equation}
 \mathrm{d}\mathrm{R}_{g_N}
 =2\Ric_{g_N}(\nabla f_N,\cdot),
 \qquad
 \Delta_{f_N}\mathrm{R}_{g_N}
 =\mathrm{R}_{g_N}-2|\Ric_{g_N}|^2.
 \label{eq:shrinker-scalar-identities}
\end{equation}
These identities follow from the contracted Bianchi identity and the soliton
equation; see \cite[Lemma~2.1]{CaoZhou2010} and the discussion following
\cite[Lemma~2.8]{LiZhang2023}.

\subsection{The product kernel and the density slice}
\label{sec:kerspace}
For the product of a compact normalized Ricci shrinker with the Gaussian factor, set
\[
 \mathcal L_{\bar f}:=\Delta_{\bar f}+2\Rm_{\bar g},\qquad
 \mathcal L_N:=\Delta_{f_N}+2\Rm_{g_N}.
\]
The eigenspaces on $N$ in the weighted divergence-free slice are
\[
 \mathscr E_\lambda(N):=
 \{\eta\in C^\infty(S^2T^*N):\delta_{X_N}\eta=0,\ \mathcal L_N\eta=\lambda\eta\}.
\]
For
$\N=\{0,1,2,\ldots\}$, let
\[
 \mathcal H_j:=\ker_{L^2(\mathrm d\gamma_d)}(\Delta_{f_E}+\frac{j}{2}),\qquad j\in\N.
\]
Thus $\mathcal H_j$ has Hermite degree $j$ and eigenvalue $-j/2$.
We set $\mathcal H_r=\{0\}$ for $r<0$.

\begin{lemma}
\label{lem:hermite-identities}
For nonnegative integers $r,s$,
\[
 \partial_{x_i}\mathcal H_r\subset\mathcal H_{r-1},\qquad
 \left(\partial_{x_i}-\frac{x_i}{2}\right)\mathcal H_r
                       \subset\mathcal H_{r+1},\qquad
 \mathcal H_r\mathcal H_s\subset\bigoplus_{k=0}^{r+s}\mathcal H_k.
\]
If $\pi_r$ denotes the orthogonal projection onto $\mathcal H_r$ in
$L^2(\mathrm d\gamma_d)$, then
\[
 \pi_r(|\vec x|^2u)=(2d+4r)u,\qquad u\in\mathcal H_r.
\]
\end{lemma}
\begin{proof}
Write $\Delta_{f_E}=\sum_i\partial_{x_i}^2-(1/2)\sum_i x_i\partial_{x_i}$.
Its commutators with the lowering and raising operators are
\[
 [\Delta_{f_E},\partial_{x_i}]=\frac{1}{2}\partial_{x_i},\qquad
 \left[\Delta_{f_E},\partial_{x_i}-\frac{x_i}{2}\right]
       =-\frac{1}{2}\left(\partial_{x_i}-\frac{x_i}{2}\right).
\]
The first two inclusions follow from the eigenvalue
$-r/2$ of $\mathcal H_r$. Hermite polynomials of degrees at
most $r+s$ span all polynomials of ordinary degree at most $r+s$,
which proves the product inclusion.
The adjoint in $L^2(\mathrm d\gamma_d)$ is
$\partial_{x_i}^*=-\partial_{x_i}+x_i/2$, so multiplication by $x_i$
is $2(\partial_{x_i}+\partial_{x_i}^*)$.
Consequently multiplication by $|\vec x|^2$ has degree-preserving part
\[
 4\sum_i\bigl(\partial_{x_i}\partial_{x_i}^*
              +\partial_{x_i}^*\partial_{x_i}\bigr)
 =-8\Delta_{f_E}+2d.
\]
Applying this identity to $u\in\mathcal H_r$ proves the last assertion.
\end{proof}

The space of infinitesimal solitonic deformations of $\bar M$ is
\[
 \ISD_{\bar M}:=
 \{h\in\operatorname{Dom}_{L^2}(\mathcal L_{\bar f}):
   \delta_{\bar X}h=0,\ \mathcal L_{\bar f}h=0\}.
\]
The operator has its self-adjoint weighted realization and the slice
condition is distributional. Elliptic regularity makes every kernel
tensor smooth.
We write $h_{NN}$ for the component of a tensor $h$ with both indices
tangent to $N$.

\begin{proposition}
\label{prop:product-kernel}
We have the finite orthogonal decomposition
\begin{equation}
 \ISD_{\bar M}=K=K_N\oplus K_+,\qquad
 K_j:=\mathcal H_j\otimes\mathscr E_{j/2}(N),\qquad
 K_N=K_0,\qquad K_+:=\bigoplus_{j\ge1}K_j.
 \label{eq:product-spectrum}
\end{equation}
Every element of $K$ has only $NN$ components.
\end{proposition}
\begin{proof}
The parallel splitting $T\bar M=TN\oplus T\R^d$ is preserved by
$\mathcal L_{\bar f}$. On the mixed and Euclidean--Euclidean blocks,
the curvature endomorphism vanishes. A zero mode in either block
therefore satisfies $\int_{\bar M}|\bar\nabla h|^2\;\mathrm d\bar\mu=0$.
The integration by parts follows by approximation in the self-adjoint
graph domain. A nonzero mixed parallel tensor would yield a nonzero
parallel vector field $Y$ on $N$, with $\Ric_{g_N}(Y,Y)=0$.
At a maximum of $f_N$, the shrinker equation gives
$\Ric_{g_N}=g_N/2-\nabla^2f_N\ge g_N/2$, a contradiction.
A Euclidean--Euclidean parallel tensor is a constant symmetric matrix
$B$; its weighted divergence is $-\frac{B\vec x}{2}$. The $NN$ block
contributes no Euclidean component to this divergence, so $B=0$.

In the remaining block, expand in Hermite modes. For
$u\in \mathcal H_j$ and a tensor $\eta$ on $N$,
\[
 \mathcal L_{\bar f}(u\eta)
 =-\frac{j}{2}u\eta+u\mathcal L_N\eta.
\]
The kernel and slice equations give $\mathcal L_N\eta=\frac{j\eta}{2}$ and
$\delta_{X_N}\eta=0$ separately at each level. The spectrum of $\mathcal L_N$ is
discrete and bounded above, with eigenvalues tending to $-\infty$.
Only finitely many values $j/2$ occur, proving the assertion.
\end{proof}

Write $J$ for the largest nonzero Hermite degree in $K$. It is fixed
by the model. The shrinker identities imply
\[
 K_G:=\{2u\Ric_{g_N}:u\in \mathcal H_2\}\subset K_2\subset K_+.
\]
Since $K_G\ne0$, we have $J\ge2$.

Set
\[
 u_*:=|\vec x|^2-2d,\qquad v_*:=2u_*\Ric_{g_N}\in K_G.
\]

\begin{definition}[Slice and inner product]
\label{def:affine-density-slice}
For $H=(h,\sigma)$ with $\bar g+h>0$ and $1+\sigma>0$, define the
metric--density chart
\[
 \mathcal C(H)=\left(\bar g+h,
 \bar f+\frac{1}{2}\log\det(I+\bar g^{-1}h)-\log(1+\sigma)\right).
\]
Thus, for $X=\mathcal C(H)$,
$\mathrm d\mu_X=(1+\sigma)\;\mathrm d\bar\mu$.
The scalar coordinate $H_{\mathrm s}=\sigma$ is a density perturbation,
whereas the scalar component of a variation of $X$ varies its potential.
The weighted divergence-free slice with normalized and centered density
is parametrized by
\[
 \mathcal S=\left\{H:\delta_{\bar X}H_{\mathrm t}=0,\
       \int_{\bar M} H_{\mathrm s}\;\mathrm d\bar\mu=0,\
       \int_{\bar M}\vec{x}H_{\mathrm s}\;\mathrm d\bar\mu=0\right\}.
\]
The scalar conditions say exactly that $\mathcal V(X)=1$ and
$\int_{\bar M}\vec{x}\;\mathrm d\mu_X=0$; the latter fixes translations.
Use the fixed $L^2(\mathrm d\bar\mu)$ inner product
\[
 \langle H,V\rangle=\int_{\bar M}
  \bigl(\langle H_{\mathrm t},V_{\mathrm t}\rangle_{\bar g}+H_{\mathrm s}V_{\mathrm s}\bigr)\;\mathrm d\bar\mu,
 \qquad
 \iota(a):=(a,0),\qquad Z:=\mathcal S\cap\iota(K)^\perp.
\]
The normal space $Z$ consists of the directions in $\mathcal S$
orthogonal to the kernel $\iota(K)$.
Let $\Pi_{\rm sl}$ be the $L^2(\mathrm d\bar\mu)$-orthogonal projection
onto tensors satisfying $\delta_{\bar X}h=0$, and set
\[
 \Pi_{S_0}\sigma=\sigma-\int_{\bar M}\sigma\;\mathrm d\bar\mu
       -\frac{1}{2}\sum_i x_i\int_{\bar M} x_i\sigma\;\mathrm d\bar\mu.
\]
Then $\Pi_{\mathcal S}H=(\Pi_{\rm sl}H_{\mathrm t},\Pi_{S_0}H_{\mathrm s})$ and
$\Pi_ZH=\Pi_{\mathcal S}H-\iota(\Pi_KH_{\mathrm t})$.
\end{definition}

\begin{lemma}
\label{lem:density-chart}
For $X=(g,f)=\mathcal C(H)$ and coordinate directions
$U=(u,\xi)$, $V=(v,\eta)$, one has
\begin{align*}
 \mathrm D\mathcal C_H[V]
   &=\left(v,\frac{1}{2}\Tr_gv-\frac{\eta}{1+\sigma}\right),
 &q_X(\mathrm D\mathcal C_H[V])&=\frac{\eta}{1+\sigma},\\
 \mathrm D^2\mathcal C_H[U,V]
   &=\left(0,-\frac{1}{2}\langle u,v\rangle_g
                         +\frac{\xi\eta}{(1+\sigma)^2}\right).
\end{align*}
In particular,
$\mathrm D\mathcal C_0[\iota(a)]=(a,\frac{\Tr_{\bar g}a}{2})$.
For each fixed integer $j\ge0$, take a sufficiently small bounded
$C^j$ neighborhood of $H=0$. For $X_i=\mathcal C(H_i)$ in this neighborhood,
\begin{equation}
 C_j^{-1}\|H_1-H_0\|_{H^j}
 \le\|X_1-X_0\|_{H^j}\le C_j\|H_1-H_0\|_{H^j}.
 \label{eq:density-physical-comparison}
\end{equation}
The same comparison holds for bounded $C^j$ norms and
pointwise derivative sums through order $j$. If
$\|H_i\|_{C^j}\le\epsilon$, then
\[
 \|X_1-X_0-\mathrm D\mathcal C_0[H_1-H_0]\|_{H^j}
       \le C_j\epsilon\|H_1-H_0\|_{H^j}.
\]
The constants are independent of support and comparison radii.
\end{lemma}
\begin{proof}
The derivative identities follow by differentiating the chart.
Integrating $\mathrm D\mathcal C-\mathrm D\mathcal C_0$ along the
segment from $H_0$ to $H_1$ gives the last estimate by Leibniz's rule,
since its coefficients are $O(\epsilon)$ in $C^j$. The same argument
applies to bounded $C^j$ norms and pointwise derivative sums.
The fixed triangular map
$\mathrm D\mathcal C_0(h,\sigma)=(h,\frac12\Tr_{\bar g}h-\sigma)$
is invertible, so absorption of the remainder proves both sides of
\eqref{eq:density-physical-comparison} and their stated analogues.
\end{proof}

For a variation $P$ of $X=(g,f)$, the entropy first variation is
\begin{equation}
 \mathrm D\mathcal W_X[P]
 =\int_{\bar M}\bigl(\langle\Phi(X),P_{\mathrm t}\rangle_g
               +(w(X)-1)q_X(P)\bigr)\;\mathrm d\mu_X;
 \label{eq:physical-entropy-variation}
\end{equation}
see \cite[Sections~1.1 and~3.1]{Perelman2002}.

\begin{definition}
\label{def:constrained-entropy-gradient}
Set
$\mathscr W=\mathcal W\circ\mathcal C$. Its constrained gradient
$\mathcal G(H)\in\mathcal S$ is characterized by
\[
 \langle\mathcal G(H),V\rangle=\mathrm D\mathscr W_H[V]
 \qquad\text{for every }V\in\mathcal S.
\]
For $H\in\mathcal S$, this is the gradient of the restricted entropy
with respect to the fixed $L^2(\mathrm d\bar\mu)$ inner product. The same identity defines
$\mathcal G(H)$ for nearby $H$ outside $\mathcal S$ by projecting the
ambient gradient onto $\mathcal S$.
The normal equation is $\Pi_Z\mathcal G(H)=0$.
\end{definition}

\begin{lemma}
\label{lem:constrained-entropy-gradient}
For $X=(g,f)=\mathcal C(H)$ and $\mathrm A=g^{-1}\bar g$, one has
\begin{equation}
 \mathcal G(H)=\left(\Pi_{\rm sl}[(1+\sigma)\mathrm A\Phi(X)\mathrm A],
                         \Pi_{S_0}(w(X)-1)\right).
 \label{eq:physical-constrained-gradient}
\end{equation}
For a general coordinate direction $H=(h,\sigma)$,
\[
 \mathrm D\mathcal G(0)[H]=
 \left(\frac{1}{2}\Pi_{\rm sl}\mathcal L_{\bar f}h,
 \Pi_{S_0}\bigl[-(2\Delta_{\bar f}+1)\sigma
       +\delta_{\bar X}\delta_{\bar X}h\bigr]\right).
\]
Its restriction to $\mathcal S$ is
\begin{equation}
 \mathfrak B_0(h,\sigma)
   =\left(\frac{1}{2}\mathcal L_{\bar f}h,
                          -(2\Delta_{\bar f}+1)\sigma\right).
 \label{eq:physical-constrained-hessian}
\end{equation}
\end{lemma}
\begin{proof}
The first variation and Lemma~\ref{lem:density-chart} give
\[
 \mathrm D\mathscr W_H[V]
 =\int_{\bar M}\bigl((1+\sigma)\langle\Phi(X),V_{\mathrm t}\rangle_g
                         +(w(X)-1)V_{\mathrm s}\bigr)\;\mathrm d\bar\mu.
\]
In a $\bar g$-orthonormal frame,
$\langle\Phi,V_{\mathrm t}\rangle_g=\langle\mathrm A\Phi\mathrm A,V_{\mathrm t}\rangle_{\bar g}$,
which proves \eqref{eq:physical-constrained-gradient}.
Lemma~\ref{lem:pair-first-variation} and the chart derivative give
\begin{align*}
 \mathrm D(\Phi\circ\mathcal C)_0[H]
  &=\frac{1}{2}\mathcal L_{\bar f}h
    +\delta_{\bar X}^*\delta_{\bar X}h+\bar\nabla^2\sigma,\\
 \mathrm D(w\circ\mathcal C)_0[H]
  &=-(2\Delta_{\bar f}+1)\sigma
             +\delta_{\bar X}\delta_{\bar X}h.
\end{align*}
Since $\Phi(\bar X)=0$ and both projections are fixed, differentiating
\eqref{eq:physical-constrained-gradient} gives the claimed derivative:
$\Pi_{\rm sl}$ annihilates the symmetric-gradient terms. On $\mathcal S$,
$\delta_{\bar X}h=0$, and the two diagonal operators preserve the
slice and scalar constraints, yielding
\eqref{eq:physical-constrained-hessian}.
\end{proof}

The bound $\lambda_1(-\Delta_{f_N})>1/2$ on $N$ and the Hermite
spectrum show that the scalar kernel consists precisely of linear
functions of the Euclidean coordinates. The first-moment condition removes these, so
$\ker\mathfrak B_0=\iota(K)$. Lemma~\ref{lem:weighted-linear}
gives the fixed projection bounds and the normal inverse.

\begin{lemma}
\label{lem:weighted-linear}
For every integer $j\ge0$, the fixed projections
$\Pi_{\mathcal S},\Pi_K,\Pi_Z$ are bounded on $H^j$, and
\[
 \mathfrak B_0:Z\cap H^{j+2}\longrightarrow Z\cap H^j
\]
is an isomorphism satisfying
\begin{equation}
 \|U\|_{H^{j+2}}\le C_j\|\mathfrak B_0U\|_{H^j}
 \qquad(U\in Z\cap H^{j+2}).
 \label{eq:full-normal-inverse}
\end{equation}
Identify Killing fields of $\bar g$ with their $\bar g$-dual one-forms.
On their orthogonal complement in $L^2(\mathrm d\bar\mu)$, the operator
$\bar{\mathcal P}=\delta_{\bar X}\delta_{\bar X}^*$ is invertible, with
the estimate, for every $\eta\in H^j$ in this complement,
\[
 \|\bar{\mathcal P}^{-1}\eta\|_{H^{j+2}}
       \le C_j\|\eta\|_{H^j}.
\]
The constants depend only on the fixed product and $j$.
\end{lemma}
\begin{proof}
The metric slice projection is bounded on every $H^j$ by
\cite[Lemma~3.17]{LiZhang2023}. The other projections add or remove
fixed finite-dimensional spaces of smooth tensors or Hermite
polynomials, and hence have the same boundedness. The tensor estimate
on the normal space follows from \cite[Theorem~3.19]{LiZhang2023}.
The scalar estimate follows from \cite[Propositions~3.7--3.8]{LiZhang2023}
and the spectral gap of $2\Delta_{\bar f}+1$ after removing the degree-one
Hermite modes. Together with \eqref{eq:physical-constrained-hessian},
these give \eqref{eq:full-normal-inverse}. The separated operators on $N$
and $\R^d$ have discrete spectra, and the kernel has been removed on $Z$;
self-adjointness therefore also gives surjectivity. The gauge inverse
and its estimate follow from \cite[Proposition~3.10 and its proof]{LiZhang2023},
whose spectral argument applies to every right-hand side orthogonal
to the one-forms dual to Killing fields.
\end{proof}

\begin{lemma}
\label{lem:bounded-gradient-comparison}
For every integer $j\ge0$, on a sufficiently small bounded
$C^{j+2}$ neighborhood of zero in metric--density coordinates,
\begin{align}
 \|(\mathrm D\mathcal G(H)-\mathrm D\mathcal G(0))U\|_{H^j}
 &\le C_j\epsilon_j\|U\|_{H^{j+2}},
 \label{eq:gradient-coefficient-comparison}\\
 \|\mathcal G(H)\|_{H^j}
 &\le C_j\|\Phi(\mathcal C(H))\|_{H^j}.
 \label{eq:gradient-by-defect}
\end{align}
Here $\|H\|_{C^{j+2}}\le\epsilon_j$.
Estimate~\eqref{eq:gradient-coefficient-comparison} holds for arbitrary
coordinate pairs $H$ and directions $U$;
estimate~\eqref{eq:gradient-by-defect} assumes
$\int_{\bar M} H_{\mathrm s}\;\mathrm d\bar\mu=0$.
\end{lemma}
\begin{proof}
Fix $j\ge0$. Differentiate \eqref{eq:physical-constrained-gradient}
and the formulas for $\mathrm D\mathcal C$ and $\mathrm D^2\mathcal C$.
Before the fixed projections, the difference from the reference
ambient differential is a sum of terms
\[
 a_2\bar\nabla^2U+a_1\bar\nabla U+a_0U,
 \qquad |\bar\nabla^i a_k|
 \le C_j\epsilon_j\langle\vec{x}\rangle^{2-k}
 \quad(0\le i\le j).
\]
Only $\mathrm df_E=\sum_i\frac{x_i}{2}\,\mathrm dx_i$ grows.
The weighted multiplication estimate
$\|\langle\vec x\rangle T\|_{H^q}\le C_q\|T\|_{H^{q+1}}$,
applied twice when necessary, and Leibniz's rule give
\eqref{eq:gradient-coefficient-comparison}.
The fixed projection bounds are supplied by
Lemma~\ref{lem:weighted-linear}; the scalar projections differ from the
identity by fixed finite-rank projections onto Hermite polynomials.
The metric--density chart comparison controls
all bounded coefficient multipliers in terms of the indicated
$C^{j+2}$ norm.

For the second estimate put $X=(g,f)=\mathcal C(H)$, $\sigma=H_{\mathrm s}$,
and
\[
 r_X=w(X)-\int_{\bar M} w(X)\;\mathrm d\bar\mu.
\]
The reference scalar spectral gap and weighted elliptic estimate give
a solution
\[
 \Delta_{\bar f}v=r_X,\qquad
 \int_{\bar M} v\;\mathrm d\bar\mu=0,\qquad
 \|v\|_{H^2}\le C\|r_X\|_{L^2}.
\]
Set $V=(1+\sigma)^{-1}\bar\nabla v$. The bounded $C^2$ neighborhood gives
\[
 \|L_Vg\|_{L^2(\mu_X)}\le C\|v\|_{H^2}\le C\|r_X\|_{L^2}.
\]
The contracted Bianchi identity $\mathrm dr_X=-2\delta_X\Phi(X)$ and
$\mathrm d\mu_X=(1+\sigma)\mathrm d\bar\mu$ give
\[
 \|r_X\|_{L^2}^2
 =-\int_{\bar M}\mathrm dr_X(\bar\nabla v)\;\mathrm d\bar\mu
 =-\int_{\bar M}\langle\Phi(X),L_Vg\rangle_g\;\mathrm d\mu_X.
\]
Cutoff approximation in $H^2$ justifies the integrations, with the weighted
multiplication estimate controlling the drift terms. Cauchy--Schwarz and
equivalence of the two weighted norms prove
$\|r_X\|_{L^2}\le C\|\Phi(X)\|_{L^2}$.

For each $j\ge1$, the formula for $\delta_X$ has bounded derivative
coefficients and drift coefficients of at most linear growth.
Leibniz's rule and the weighted multiplication estimate thus give
\[
 \|\delta_X\Phi(X)\|_{H^{j-1}}
       \le C_j\|\Phi(X)\|_{H^j}.
\]
Using the first-order Bianchi identity directly, we obtain
\[
 \|r_X\|_{H^j}
 \le C_j\bigl(\|r_X\|_{L^2}+\|\mathrm dr_X\|_{H^{j-1}}\bigr)
 \le C_j\|\Phi(X)\|_{H^j}.
\]
Since
$\Pi_{S_0}(w-1)=\Pi_{S_0}r_X$, the fixed projection bound controls
the scalar entry of $\mathcal G$. The bounded coefficient multipliers
$(1+H_{\mathrm s})\mathrm A(\cdot)\mathrm A$ and the fixed slice
projection control its tensor entry in
\eqref{eq:physical-constrained-gradient}, proving
\eqref{eq:gradient-by-defect}.
\end{proof}

\section{Normal reduction for generalized cylinders}
\label{sec:quantitative}
\subsection{Reduction on the compact factor}
\label{sec:local-kuranishi}
Use the metric--density chart on $N$ and the spaces
\[
 \mathcal S_N=\left\{H:\delta_{X_N}H_{\mathrm t}=0,\
                    \int_NH_{\mathrm s}\;\mathrm d\mu_N=0\right\},
 \qquad Z_N=\mathcal S_N\cap\iota(K_N)^\perp.
\]
Here the inner product, chart, and $\iota$ are defined on $N$ as above.
Write $\mathcal G_N$ for the constrained gradient
\eqref{eq:physical-constrained-gradient}, with the scalar projection
removing only the constant mode.

\begin{lemma}
\label{lem:compact-variational-graph}
For $b\in K_N$ sufficiently small, there is a smooth family of pairs
$v_0(b)$ on $N$, depending real analytically on $b$ in each fixed
$C^q$ topology, such that
\begin{equation}
 \begin{gathered}
 H_b=\iota(b)+v_0(b),\qquad v_0(b)\in Z_N,\qquad
 \Pi_{Z_N}\mathcal G_N(H_b)=0,\\
 v_0(0)=\mathrm Dv_0(0)=0.
 \end{gathered}
 \label{eq:variational-compact-graph}
\end{equation}
This is the unique nearby solution with kernel coordinate $b$, and
$\|v_0(b)\|_{C^q(N)}\le C_q|b|^2$ for every integer $q\ge0$.
\end{lemma}
\begin{proof}
Fix an integer $q\ge1$ and set
\[
 F(b,v)=\Pi_{Z_N}\mathcal G_N(\iota(b)+v).
\]
The slice projection on $N$ is bounded on H\"older spaces by
\cite[Proposition~2.3]{LiZhang2023}. The Hessian
\eqref{eq:physical-constrained-hessian} therefore makes
$\mathrm D_vF(0,0)=\mathfrak B_0|_{Z_N}$ an isomorphism from
$Z_N\cap C^{q+2,1/2}$ to $Z_N\cap C^{q,1/2}$:
the tensor kernel has been removed, and
$\lambda_1(-\Delta_{f_N})>1/2$ excludes a mean-zero scalar kernel.
The map $F$ is analytic and $F(0,0)=0$, so the analytic implicit
function theorem gives the unique nearby solution $v_0(b)$.
Since $\mathrm D_bF(0,0)=0$, differentiation of the normal equation
gives $\mathrm Dv_0(0)=0$, and hence
$\|v_0(b)\|_{C^{q+2,1/2}}\le C_q|b|^2$.
Elliptic regularity and uniqueness identify these solutions at all
regularity orders. The continuous embeddings into $C^q$ give the
stated analyticity and estimates for every integer $q\ge0$.
\end{proof}

Extend $H_b$ constantly in $\vec{x}$ and put
\begin{equation}
 X_b:=\mathcal C(H_b),\qquad
 e(b):=\mathscr W(H_b)=\mathcal W(X_b).
 \label{eq:Xb-background}
\end{equation}
The Gaussian factor contributes zero to the normalized entropy.
The product extension has mass one and zero first moments. Since
$v_0(b)\perp\iota(K_N)$ and the chart is affine in its tensor
component, $\Pi_K(X_b-\bar X)_{\mathrm t}=b$.

\begin{proposition}
\label{prop:compact-transfer}
The potential of $X_b$ on $N$ is the unique normalized minimizing potential realizing
$\boldsymbol{\mu}((X_b)_{\mathrm t}|_{TN},1)$. For every integer $q\ge0$,
\begin{equation}
 \|\Phi(X_b)\|_{C^q(N)}
 \le C_q|\mathrm De(b)|.
 \label{eq:variational-full-defect}
\end{equation}
In particular, $X_b$ is a Ricci shrinker if and only if
$\mathrm De(b)=0$.
\end{proposition}
\begin{proof}
Test \eqref{eq:variational-compact-graph} against $(0,\eta)$ with
$\int_N\eta\;\mathrm d\mu_N=0$. Such variations lie in $Z_N$, so
\eqref{eq:physical-constrained-gradient} gives
\[
 \int_N(w(X_b)-1)\eta\;\mathrm d\mu_N=0.
\]
Thus $w(X_b)$ is spatially constant, and integration by parts gives
$w(X_b)=e(b)$. By \cite[Lemma~2.2 and its proof]{SunWang2015},
the nearby normalized solution of this scalar Euler--Lagrange equation is the
unique minimizing potential. This proves the first assertion.

Fix $q\ge1$. Write $\Psi_b=\Phi(X_b)|_{TN}$,
$\vartheta_b=1+(H_b)_{\mathrm s}$, and
$\mathrm A_b=((X_b)_{\mathrm t}|_{TN})^{-1}g_N$. The normal equation implies
$\mathcal G_N(H_b)\in\iota(K_N)$. Differentiating $e(b)$ pairs
this gradient with
$\iota(\beta)+\mathrm Dv_0(b)[\beta]$; the second summand is
normal and contributes zero. Hence
$\mathcal G_N(H_b)=\iota(\nabla e(b))$, where $\nabla e$ uses
the fixed $L^2(\mathrm d\mu_N)$ inner product on $K_N$, and
\begin{equation}
 \Pi_{\rm sl}(\vartheta_b \mathrm A_b\Psi_b \mathrm A_b)=\nabla e(b).
 \label{eq:compact-variational-gradient}
\end{equation}
The contracted Bianchi identity
$\mathrm dw(X_b)=-2\delta_{X_b}\Phi(X_b)$ gives
$\delta_{X_b}\Phi(X_b)=0$. Restricting to $N$, use the Hodge
decomposition for $(g_N,f_N)$ and its inverse on the
$L^2(\mathrm d\mu_N)$-orthogonal complement of the $g_N$-dual one-forms
of Killing fields. These give
\[
 \|(I-\Pi_{\rm sl})\Psi_b\|_{C^{q,1/2}}
 \le C\|\delta_{X_N}\Psi_b\|_{C^{q-1,1/2}}
 \le C|b|\|\Psi_b\|_{C^{q,1/2}}.
\]
The multiplier $\vartheta_b \mathrm A_b(\cdot)\mathrm A_b-I$ is $O(|b|)$ in every
fixed coefficient norm. By \eqref{eq:compact-variational-gradient},
\[
 \|\Pi_{\rm sl}\Psi_b-\nabla e(b)\|_{C^{q,1/2}}
 \le C|b|\|\Psi_b\|_{C^{q,1/2}}.
\]
Finite-dimensional norm equivalence and absorption give the
$C^{q,1/2}$ bound. The continuous embedding into $C^q$ gives
\eqref{eq:variational-full-defect}, including $q=0$. This estimate and
\eqref{eq:compact-variational-gradient} give the last assertion.
\end{proof}

\begin{assumption}[Finite-order rigidity of compact shrinkers]
\label{ass:compact-rigidity}
There exist an integer $\ell\ge1$ and constants $c_K,\epsilon_K>0$
such that
\begin{equation}
 |\mathrm De(b)|\ge c_K|b|^\ell
 \qquad\text{whenever }|b|\le\epsilon_K.
 \label{eq:compact-rigidity}
\end{equation}
The convention $\ell=1$ means $K_N=0$; if $K_N\ne0$, we take
$\ell\ge2$.
\end{assumption}

\begin{proposition}
\label{prop:compact-isolation}
Assumption~\ref{ass:compact-rigidity} for some finite $\ell$ is
equivalent to local isolation of $X_N$ among normalized compact
Ricci shrinkers, modulo diffeomorphisms.
\end{proposition}
\begin{proof}
The compact slice theorem \cite{PodestaSpiro2013}, in the
fixed-divergence formulation of \cite[Theorem~2.4]{LiZhang2023},
puts every sufficiently close metric in the fixed slice. A nearby normalized shrinker, expressed in metric--density
coordinates, solves the normal equation. Uniqueness in
Lemma~\ref{lem:compact-variational-graph} places it on the family
$H_b$. By Proposition~\ref{prop:compact-transfer}, the nearby
shrinker equation is therefore equivalent to $\mathrm De(b)=0$.
The inequality \eqref{eq:compact-rigidity} makes this zero set
$\{0\}$. Conversely, local isolation makes $0$ an isolated zero of
$\mathrm De$. Applying the real analytic distance inequality
\cite[Theorem~6.4 and Remark~6.5]{BierstoneMilman1988} to
$|\mathrm De|^2$ gives $|b|^p\le C|\mathrm De(b)|$ for some finite
$p>0$. Increasing $p$ to an integer proves the assumption on a
smaller ball; when $K_N=0$ the assumption holds with $\ell=1$.
\end{proof}

\begin{remark}
This is also the compact rigidity characterization in
\cite[Theorem~1.1 and Corollaries~2.27, 2.28]{LiZhang2023}.
\end{remark}

\subsection{Radial variations along the family \texorpdfstring{$X_b$}{Xb}}
\label{sec:local-moving-gauge}

Let $H_b=\iota(b)+v_0(b)$, $v_0(b)\in Z$, denote the family
constructed in Lemma~\ref{lem:compact-variational-graph}, extended
constantly in $\vec{x}$. Its normalized pair is
$X_b=\mathcal C(H_b)$, and $e(b)=\mathcal W(X_b)$.
The construction on $N$ gives, for every fixed integer $q\ge0$,
\begin{equation}
 \mathcal G(H_b)=\iota(\nabla e(b)),\qquad
 w(X_b)=e(b),\qquad
 \|\Phi(X_b)\|_{C^q(N)}\le C_q|\mathrm De(b)|.
 \label{eq:variational-compact-defect}
\end{equation}

Write $X_b=(g_b+g_E,f_b+f_E)$ and set $\Phi_b=\Phi(X_b)$.
The tensor $\Phi_b$ has only $NN$ components and is independent of $\vec{x}$.

\begin{lemma}
\label{lem:moving-gauge-field}
For sufficiently small $b$, there is a unique vector field
\[
 V_b\in\bigl(\Gamma(TN)\otimes\mathcal H_2\bigr)
       \oplus\bigl(C^\infty(N;\mathbb R^d)\otimes\mathcal H_1\bigr),
 \qquad V_b^{\flat_{\bar g}}\perp\ker\bar{\mathcal P},
\]
satisfying
\begin{equation}
 \delta_{\bar X}L_{V_b}(g_b+g_E)
       =-\delta_{\bar X}(2u_*\Ric_{g_b}).
 \label{eq:moving-gauge-equation}
\end{equation}
Its coefficients on $N$ depend real analytically on $b$ in every fixed
$C^q$ topology, and $V_0=0$. For every integer $q\ge0$ and every $A\ge0$,
\begin{equation}
 \sup_{\bar M}\langle\vec x\rangle^{-2}[V_b]_q
 +\sum_{i=0}^q\|\langle\vec x\rangle^A\bar\nabla^iV_b\|_{L^2}
 \le C_{q,A}|b|.
 \label{eq:moving-gauge-bounds}
\end{equation}
The orthogonality and norms use the fixed reference metric and measure.
\end{lemma}
\begin{proof}
Let $\mathcal A_bV=\delta_{\bar X}L_V(g_b+g_E)$, identifying vector
fields and one-forms by $\bar g$. If the components of $V$ tangent
to $N$ and $\mathbb R^d$ have Hermite degrees two and one, respectively,
then the $NN$, mixed, and Euclidean--Euclidean components
of $L_V(g_b+g_E)$ have degrees two, one, and zero. Differentiation
along $N$ preserves these degrees. The Euclidean weighted
derivative $\partial_{x_i}-\frac{x_i}{2}$ raises Hermite degree by one
by Lemma~\ref{lem:hermite-identities}, so
$\mathcal A_b$ preserves the stated space of vector fields. The right-hand
side of \eqref{eq:moving-gauge-equation} lies in the same space.
For every $\omega\in\ker\bar{\mathcal P}$, adjointness gives
\[
 \langle\mathcal A_bV,\omega\rangle
 =\langle L_V(g_b+g_E),\delta_{\bar X}^*\omega\rangle=0;
\]
the right-hand side is orthogonal to these forms for the same reason.

At $b=0$,
\[
 \mathcal A_0V=-2\bar{\mathcal P}(V^{\flat_{\bar g}}).
\]
Lemma~\ref{lem:weighted-linear} gives its inverse on the intersection
of this space with the weighted $L^2$-orthogonal complement of
$\ker\bar{\mathcal P}$. The only Killing fields of $\bar g$ in this
space are Euclidean rotations: Killing fields on $N$ and Euclidean
translations have degree zero. The finite Hermite
expansion identifies its Sobolev norms with norms of finitely many
coefficients on $N$. Since $g_b-g_N=O(|b|)$ in every fixed $C^q$
norm, for every fixed $q$ one has
\begin{align*}
 \|(\mathcal A_b-\mathcal A_0)V\|_{H^q}
     &\le C_q|b|\|V\|_{H^{q+2}},\\
 \|\delta_{\bar X}(2u_*\Ric_{g_b})\|_{H^q}
     &\le C_q|b|.
\end{align*}
For the second estimate, the shrinker identity
$\mathrm d\mathrm R_{g_N}=2\Ric_{g_N}(\nabla f_N,\cdot)$ and the
contracted Bianchi identity give
\[
 \delta_{X_N}\Ric_{g_N}
 =\frac{1}{2}\mathrm d\mathrm R_{g_N}
       -\Ric_{g_N}(\nabla f_N,\cdot)=0.
\]
Since $u_*$ depends only on $\vec x$ and $\Ric_{g_N}$ has only $NN$ components,
$\delta_{\bar X}(2u_*\Ric_{g_N})=2u_*\delta_{X_N}\Ric_{g_N}=0$.
Thus the right-hand side vanishes at $b=0$.
A Neumann series on the same orthogonal complement gives a unique
solution with $\|V_b\|_{H^{q+2}}\le C_q|b|$.

In the fixed Hermite basis the equation is a uniformly elliptic
system for the coefficients on $N$. Elliptic regularity and uniqueness
identify the solutions at different orders and give every coefficient
$C^q$ bound. The operator and right-hand side depend analytically on $b$,
so the inverse does as well in each fixed topology. Each Hermite
factor has degree at most two. Its derivatives are bounded by
$C_q\langle\vec x\rangle^2$, and all its polynomially weighted Gaussian
moments are finite. These observations prove
\eqref{eq:moving-gauge-bounds}. The vanishing of the right-hand side at $b=0$
and uniqueness give $V_0=0$.
\end{proof}

\begin{proposition}
\label{prop:moving-radial-hessian}
Let $V_b$ be the field in Lemma~\ref{lem:moving-gauge-field}, and set
$\vartheta_b=1+(H_b)_{\mathrm s}$. Define
\begin{align*}
 J_b^0&=(2u_*\Ric_{g_b},u_*\mathrm R_{g_b}),\\
 J_b&=J_b^0+
       \bigl(L_{V_b}(g_b+g_E),V_b(f_b+f_E)\bigr),\\
 j_b&=\bigl((J_b)_{\mathrm t},\vartheta_bq_{X_b}(J_b)\bigr).
\end{align*}
Then
\[
 \mathrm D\mathcal C_{H_b}[j_b]=J_b,\qquad
 j_b\in\mathcal S\cap\iota(K_N)^\perp,\qquad
 j_0=\iota(v_*).
\]
Set $\mathfrak B_b=\mathrm D\mathcal G(H_b)|_{\mathcal S}$.
For every integer $q\ge0$, on a sufficiently small parameter neighborhood,
\begin{equation}
 \|j_b-j_0\|_{H^q}\le C_q|b|,\qquad
 \|\mathfrak B_bj_b\|_{H^q}\le C_q|\mathrm De(b)|.
 \label{eq:moving-radial-full}
\end{equation}
\end{proposition}
\begin{proof}
We first record the calculation on $N$. In this calculation all
geometric quantities use $g_b$, and we write $f=f_b$. Taking the
trace of $\Phi_b=\frac{g_b}{2}-\Ric-\nabla^2f$ gives
\[
 \Tr\Phi_b=\frac{m}{2}-\mathrm R-\Delta f.
\]
The scalar equation in \eqref{eq:variational-compact-defect} is constant
on $N$. Differentiating it yields
\[
 2\,\mathrm d\Delta f
 -2\nabla^2f(\nabla f,\cdot)+\mathrm d\mathrm R+\mathrm df=0.
\]
Consequently,
\begin{align*}
 -\mathrm d\Tr\Phi_b
 &=\mathrm d\mathrm R+\mathrm d\Delta f
   =\frac{1}{2}\mathrm d\mathrm R+
     \nabla^2f(\nabla f,\cdot)-\frac{1}{2}\mathrm df,\\
 \Phi_b(\nabla f,\cdot)
 &=\frac{1}{2}\mathrm df-\Ric(\nabla f,\cdot)
                          -\nabla^2f(\nabla f,\cdot).
\end{align*}
Adding these equalities and applying the contracted Bianchi identity gives
\begin{equation}
 \delta_{g_b,f_b}\Ric_{g_b}
 =-\mathrm d\Tr_{g_b}\Phi_b+\Phi_b(\nabla f_b,\cdot).
 \label{eq:moving-Ricci-divergence}
\end{equation}

The Ricci variation satisfies
\begin{gather*}
 \mathrm D\Ric[g_b]=0,\qquad
 \mathrm D\Ric[L_Yg_b]=L_Y\Ric,\\
 -2\mathrm D\Ric[\Ric]
       =\Delta\Ric+2\Rm(\Ric)-2\Ric\circ\Ric.
\end{gather*}
Apply these identities to $2\Phi_b=g_b-2\Ric-L_{\nabla f}g_b$.
Since
\[
 L_{\nabla f}\Ric=\nabla_{\nabla f}\Ric
      +\nabla^2f\circ\Ric+\Ric\circ\nabla^2f,
\]
we obtain
\begin{align*}
 2\mathrm D\Ric[\Phi_b]
 &=\mathcal L_f\Ric-2\Ric\circ\Ric
       -\nabla^2f\circ\Ric-\Ric\circ\nabla^2f\\
 &=\mathcal L_f\Ric-\Ric+\Phi_b\circ\Ric+\Ric\circ\Phi_b.
\end{align*}
Thus
\begin{equation}
 \mathcal L_{f_b}\Ric_{g_b}-\Ric_{g_b}
 +\Phi_b\circ\Ric_{g_b}+\Ric_{g_b}\circ\Phi_b
       =2\mathrm D\Ric_{g_b}[\Phi_b].
 \label{eq:moving-Ricci-eigenvalue}
\end{equation}

For the uncorrected variation $J_b^0$, one has
$q_{X_b}(J_b^0)=0$ and $\Delta_{f_E}u_*=-u_*$. Moreover,
\begin{align*}
 \mathcal L_{f_b+f_E}(2u_*\Ric_{g_b})
       &=2u_*(\mathcal L_{f_b}\Ric_{g_b}-\Ric_{g_b}),\\
 \delta_{X_b}(2u_*\Ric_{g_b})
       &=2u_*\delta_{g_b,f_b}\Ric_{g_b}.
\end{align*}
Lemma~\ref{lem:pair-first-variation} and
\eqref{eq:moving-Ricci-eigenvalue} therefore give
\begin{align}
 \mathrm D\Phi_{X_b}[J_b^0]
 &=2u_*\mathrm D\Ric_{g_b}[\Phi_b]
       +2\delta_{X_b}^*(u_*\delta_{g_b,f_b}\Ric_{g_b}),
       \label{eq:moving-raw-Phi}\\
 \mathrm Dw_{X_b}[J_b^0]
 &=2u_*\delta_{g_b,f_b}\delta_{g_b,f_b}\Ric_{g_b}
       +2u_*\langle\Phi_b,\Ric_{g_b}\rangle_{g_b}.
       \label{eq:moving-raw-w}
\end{align}
The Ricci derivative is a second-order operator with uniformly bounded
coefficients on $N$. Equations~\eqref{eq:moving-Ricci-divergence}
and~\eqref{eq:variational-compact-defect} control every coefficient on
$N$ in these expressions by $C_q|\mathrm De(b)|$. The remaining
factors are polynomials in $\vec x$. Hence, for every fixed $q$,
\[
 \|\mathrm D\Phi_{X_b}[J_b^0]\|_{H^q}
 +\|\mathrm Dw_{X_b}[J_b^0]\|_{H^q}
 \le C_q|\mathrm De(b)|.
\]
The same bounds hold with any fixed polynomial weight on each
covariant derivative.

The chart derivative in Lemma~\ref{lem:density-chart} gives
$\mathrm D\mathcal C_{H_b}[j_b]=J_b$. The gauge equation gives
$\delta_{\bar X}(j_b)_{\mathrm t}=0$. Also,
\[
 q_{X_b}(J_b^0)=0,\qquad
 q_{X_b}(J_b-J_b^0)=\delta_{X_b}V_b,
 \qquad (j_b)_{\mathrm s}=\vartheta_b\delta_{X_b}V_b.
\]
By Lemma~\ref{lem:hermite-identities}, Euclidean weighted divergence
raises Hermite degree by one, while differentiation along $N$ preserves it.
Thus $\delta_{X_b}V_b$ has degree two. Since $\vartheta_b$ depends
only on $N$, $(j_b)_{\mathrm s}$ has degree two as well, so its mean
and first moments vanish. Hence $j_b\in\mathcal S$.
The $NN$ component of $(j_b)_{\mathrm t}$ has degree two,
so it is orthogonal to $K_N$. Since $V_0=0$ and
$q_{\bar X}(J_0^0)=0$, one has
$j_0=(2u_*\Ric_{g_N},0)=\iota(v_*)$.
The coefficient bounds for $X_b$ and
\eqref{eq:moving-gauge-bounds} now give
$\|j_b-j_0\|_{H^q}\le C_q|b|$.

By naturality,
\[
 \mathrm D\Phi_{X_b}[J_b]
       =\mathrm D\Phi_{X_b}[J_b^0]+L_{V_b}\Phi_b,\qquad
 \mathrm Dw_{X_b}[J_b]=\mathrm Dw_{X_b}[J_b^0].
\]
The second equality uses the spatial constancy of $w(X_b)$.
The new Lie derivative satisfies
\mbox{$\|L_{V_b}\Phi_b\|_{H^q}\le C_q|b|\,|\mathrm De(b)|$}.
Differentiating \eqref{eq:physical-constrained-gradient} in the
coordinate direction $j_b$ yields the terms involving
$\mathrm D\Phi_{X_b}[J_b]$ and $\mathrm Dw_{X_b}[J_b]$ estimated above.
The remaining terms differentiate the inverse-metric or density
factors and each contains $\Phi_b$. Their coefficients on $N$ are
uniformly bounded, and their polynomial factors in $\vec x$ have fixed
degree. Estimate~\eqref{eq:variational-compact-defect}, finite Gaussian
moments, and
the projection bounds in Lemma~\ref{lem:weighted-linear} therefore give
$\|\mathfrak B_bj_b\|_{H^q}\le C_q|\mathrm De(b)|$, proving
\eqref{eq:moving-radial-full}.
\end{proof}

\subsection{Normal reduction through second order}
\label{sec:normal-reduction}

The family $H_b$ satisfies
\[
 H_b=\iota(b)+v_0(b),\qquad
 \Pi_Z\mathcal G(H_b)=0,\qquad
 \mathcal G(H_b)=\iota(\nabla e(b)).
\]
For $a\in K_+$, we construct the linear and quadratic normal corrections
to $H_b+\iota(a)$ so that
\[
 \Pi_Z\mathcal G\bigl(H_b+\iota(a)+w(b,a)\bigr)=0,
 \qquad w(b,a)\in Z,\quad w(b,0)=0,
\]
holds through quadratic order in $a$. The construction uses finite
Hermite expansions and a quadratic Taylor polynomial in $a$, with
exact dependence on $b$ and normal corrections in the infinite-dimensional
space $Z$.

\begin{definition}
\label{def:shifted-hermite-blocks}
For an integer $l\ge0$, let
$\mathcal P_l$ consist of smooth coordinate pairs $U=(h,\sigma)$ with
\[
\begin{aligned}
 h_{NN}&\in C^\infty(S^2T^*N)\otimes\mathcal H_l,\\
 h_{NE}&\in C^\infty(T^*N\otimes(\mathbb R^d)^*)
                       \otimes\mathcal H_{l-1},\\
 h_{EE}&\in C^\infty(N;S^2(\mathbb R^d)^*)\otimes\mathcal H_{l-2},\\
 \sigma&\in C^\infty(N)\otimes\mathcal H_l.
\end{aligned}
\]
Here $h_{NE}$ specifies both mixed components by symmetry. Put
\[
 \mathcal P_{\le L}=\bigoplus_{l=0}^L\mathcal P_l,
 \qquad Z_l=Z\cap\mathcal P_l.
\]
We call $\mathcal P_l$ a shifted Hermite block: the tensor degrees
$l,l-1,l-2$ decrease with the number of Euclidean indices, while the
scalar degree is $l$. Its elements have smooth coefficients on $N$ and
finite Hermite expansions. When using Sobolev domains, we take the
closures of these spaces in the indicated $H^q$ spaces, keeping the
same notation.
\end{definition}

The spaces $\mathcal P_l$ are mutually orthogonal for the fixed pairing. Choose
orthonormal bases of the finitely many Hermite factors in
$\mathcal P_{\le L}$. Gaussian integration and differentiation then
identify each weighted $H^q$ norm with an equivalent norm on the finite
collection of $H^q(N)$ spaces for its coefficients. We will use
\[
 \iota(K_l)\subset\mathcal P_l,\qquad
 \mathcal S\cap\mathcal P_0=\mathcal S_N,\qquad
 j_b\in\mathcal S\cap\mathcal P_2.
\]
The slice $\mathcal S_N$ in the middle identity is extended constantly in
$\vec x$; the last identity follows from
Lemma~\ref{lem:moving-gauge-field} and
Proposition~\ref{prop:moving-radial-hessian}.

\begin{lemma}
\label{lem:hermite-block-calculus}
The weighted divergence maps tensor components of degrees
$(l,l-1,l-2)$ to one-form components of degrees $(l,l-1)$.
Its adjoint maps these one-form degrees back to the tensor degrees;
weighted divergence of such a one-form has scalar degree $l$.
The Hessian of a scalar in $\mathcal H_l$ has tensor degrees
$(l,l-1,l-2)$.
The inverse of $\bar{\mathcal P}$ on the weighted $L^2$-orthogonal
complement of its kernel preserves these one-form spaces. Consequently
the fixed projections $\Pi_{\mathcal S},\Pi_K,\Pi_Z$ preserve
$\mathcal P_l$,
with $\Pi_K$ acting on the tensor component.
\end{lemma}
\begin{proof}
Use indices $a,b$ on $N$ and Euclidean indices $i,j$. The product
connection and measure give
\begin{align*}
 (\delta_{\bar X}h)_a
   &=(\delta_{X_N}h_{NN})_a
       +\sum_i\left(\partial_{x_i}-\frac{x_i}{2}\right)h_{ia},\\
 (\delta_{\bar X}h)_i
   &=\nabla_N^ah_{ai}-(\nabla_N^af_N)h_{ai}
       +\sum_j\left(\partial_{x_j}-\frac{x_j}{2}\right)h_{ji}.
\end{align*}
By Lemma~\ref{lem:hermite-identities}, these expressions have
degrees $l$ and $l-1$. For a one-form
$\omega=(\omega_N,\omega_E)$ of those degrees,
\begin{align*}
 (\delta_{\bar g}^*\omega)_{ab}
   &=-\frac{1}{2}(\nabla_a\omega_b+\nabla_b\omega_a),\\
 (\delta_{\bar g}^*\omega)_{ai}
   &=-\frac{1}{2}(\nabla_a\omega_i+\partial_{x_i}\omega_a),\\
 (\delta_{\bar g}^*\omega)_{ij}
   &=-\frac{1}{2}(\partial_{x_i}\omega_j+\partial_{x_j}\omega_i),\\
 \delta_{\bar X}\omega
   &=\delta_{X_N}\omega_N+
       \sum_i\left(\partial_{x_i}-\frac{x_i}{2}\right)\omega_i.
\end{align*}
The adjoint here is also $\delta_{\bar X}^*$. These formulas prove the
adjoint and scalar divergence assertions. The Hessian components of
a scalar $v$ are $\nabla_a\nabla_bv$, $\nabla_a\partial_{x_i}v$,
and $\partial_{x_i}\partial_{x_j}v$, proving the remaining degree claim.

It follows that $\bar{\mathcal P}=\delta_{\bar X}\delta_{\bar X}^*$
preserves each of these one-form spaces. Their orthogonal projections
commute with this self-adjoint operator and with the projection onto
its kernel. Uniqueness in Lemma~\ref{lem:weighted-linear} therefore
shows that its inverse on the orthogonal complement of the kernel
preserves these spaces as well. For the corresponding tensor components,
the Hodge formula
\[
 (I-\Pi_{\rm sl})h
   =\delta_{\bar X}^*\bar{\mathcal P}^{-1}\delta_{\bar X}h
\]
then proves preservation by $\Pi_{\rm sl}$. The right-hand side is
compatible with the inverse because weighted divergence is orthogonal
to $\ker\bar{\mathcal P}$. The tensor kernel splits according to
$\mathcal P_l$, so $\Pi_K$ also preserves these spaces. Finally,
$\Pi_{S_0}$ removes the constant functions at $l=0$ and the linear
functions of $\vec x$ at $l=1$. It leaves the other spaces unchanged.
This proves the assertions for
$\Pi_{\mathcal S}$ and $\Pi_Z$.
\end{proof}

\begin{lemma}
\label{lem:moving-normal-inverse}
The operator $\mathfrak B_b=\mathrm D\mathcal G(H_b)|_{\mathcal S}$
is symmetric for the fixed pairing and preserves
$\mathcal S\cap\mathcal P_l$. The projected operator
\[
 \mathfrak B_{b,Z}:=\Pi_Z\mathfrak B_b|_Z
\]
preserves $Z_l$. On a fixed sufficiently small parameter neighborhood,
\[
 \mathfrak B_{b,Z}:Z\cap H^2\longrightarrow Z\cap L^2
\]
is an isomorphism with a uniformly bounded inverse. For each fixed
finite $L$ and integer $q\ge0$, after shrinking the neighborhood
depending on $q,L$, this inverse restricts to the Sobolev
completions of $Z\cap\mathcal P_{\le L}$ and satisfies
\begin{equation}
 \|\mathfrak B_{b,Z}^{-1}V\|_{H^{q+2}}
       \le C_{q,L}\|V\|_{H^q}.
 \label{eq:moving-normal-inverse}
\end{equation}
On each such neighborhood the inverse depends analytically on $b$
as an operator from $H^q$ to $H^{q+2}$.
\end{lemma}
\begin{proof}
Put $\vartheta_b=1+(H_b)_{\mathrm s}$ and
$\mathrm A_b=(g_b+g_E)^{-1}\bar g$. For $U=(h,\sigma)$, write
$J=\mathrm D\mathcal C_{H_b}[U]$. The chart derivative gives
\[
 q_{X_b}(J)=\frac{\sigma}{\vartheta_b},\qquad
 J_{\mathrm s}=\frac{1}{2}\Tr_{g_b+g_E}h
                         -\frac{\sigma}{\vartheta_b}.
\]
If $U\in\mathcal P_l$, the first expression has degree $l$, whereas
the second can contain degrees $l$ and $l-2$. The product metric
$g_b+g_E$ and potential $f_b+f_E$ have coefficients depending only
on $N$, apart from $f_E$. Thus the component formulas in
Lemma~\ref{lem:hermite-block-calculus} apply equally to this product.
Lemma~\ref{lem:pair-first-variation} shows that
$\mathrm D\Phi_{X_b}[J]$ has tensor degrees $(l,l-1,l-2)$ and
$\mathrm Dw_{X_b}[J]$ has scalar degree $l$.

To include the differentiated metric and density factors, set
$\dot{\mathrm A}_b=-\mathrm A_bh\mathrm A_b$ in a
$\bar g$-orthonormal frame. Differentiating
\eqref{eq:physical-constrained-gradient} gives
\begin{align*}
 (\mathfrak B_bU)_{\mathrm t}
   &=\Pi_{\rm sl}\bigl[\vartheta_b\mathrm A_b
                   \mathrm D\Phi_{X_b}[J]\mathrm A_b
                    +\sigma\mathrm A_b\Phi_b\mathrm A_b
       +\vartheta_b(\dot{\mathrm A}_b\Phi_b\mathrm A_b
                         +\mathrm A_b\Phi_b\dot{\mathrm A}_b)\bigr],\\
 (\mathfrak B_bU)_{\mathrm s}
   &=\Pi_{S_0}\mathrm Dw_{X_b}[J].
\end{align*}
The factors $\vartheta_b$ and $\mathrm A_b$ depend only on $N$,
$\mathrm A_b$ is block diagonal, and $\Phi_b$ has only $NN$ components.
The terms containing $\Phi_b$ therefore have the same tensor degrees
$(l,l-1,l-2)$. The fixed projections preserve these degrees,
proving that $\mathfrak B_b$ preserves $\mathcal S\cap\mathcal P_l$
and $\mathfrak B_{b,Z}$ preserves $Z_l$.

Since $\|H_b\|_{C^{q+2}}\le C_q|b|$, Lemma~\ref{lem:bounded-gradient-comparison}
gives, for every fixed integer $q\ge0$ and sufficiently small $b$,
\begin{equation}
 \|(\mathfrak B_b-\mathfrak B_0)U\|_{H^q}
       \le C_q|b|\|U\|_{H^{q+2}},
 \qquad U\in\mathcal S\cap H^{q+2}.
 \label{eq:moving-hessian-perturbation}
\end{equation}

The operator $\mathfrak B_b$ is the entropy Hessian on $\mathcal S$
for the fixed pairing. Cutoff approximation and the coefficient bounds
underlying \eqref{eq:gradient-coefficient-comparison} extend its symmetry
to the $H^2$ domain, so the projected operator $\mathfrak B_{b,Z}$ is
symmetric as well.
Lemma~\ref{lem:weighted-linear} and
\eqref{eq:moving-hessian-perturbation} give
\[
 \| (\mathfrak B_{b,Z}-\mathfrak B_{0,Z})
                  \mathfrak B_{0,Z}^{-1}\|_{H^q\to H^q}
       \le C_q|b|.
\]
Thus $\mathfrak B_{b,Z}\mathfrak B_{0,Z}^{-1}$ is invertible by a
Neumann series for sufficiently small $|b|$, depending on $q$.
At $q=0$ this proves the global $H^2$ isomorphism. For each fixed
$q$ and finite $L$, preservation of $Z_l$ gives the restriction and
estimate \eqref{eq:moving-normal-inverse} after shrinking the
neighborhood. Uniqueness at the $H^2$ level identifies the inverses
obtained at different orders. The coefficient families and the
Neumann series are analytic in $b$ between these fixed Sobolev spaces.
\end{proof}

For a function of $(b,a)$, let $\mathrm T_a^k$ denote its Taylor polynomial
of degree at most $k$ at $a=0$, including the constant term, with
$b$ fixed. On polynomial coordinate pairs, expressions involving
$\mathcal G$ or $\mathscr W$ below mean their finite Taylor
coefficients. Thus $O_a(|a|^3)$ means that the coefficients of
degrees zero, one, and two vanish.

The coefficients are obtained by differentiating the local coordinate
expressions at $H_b$, integrating the entropy coefficients against
$\mathrm d\bar\mu$, and applying the fixed projections for $\mathcal G$.
Their polynomial growth gives Gaussian integrability. To justify these
operations, multiply the polynomial directions by
$\chi(\vec x/R)$, differentiate the ambient entropy at fixed $R$,
and then let $R\to\infty$. All differentiated cutoff errors vanish
by Gaussian tail estimates.

\begin{proposition}
\label{prop:second-order-normal-reduction}
Let $J\ge2$ be the largest integer $j$ for which $K_j\ne0$, and set
$L_1=2J+4$. For every integer $q\ge0$, on a sufficiently small
parameter neighborhood depending on $q$, the following construction
is well defined for $b\in K_N$ and $a\in K_+$.
The polynomial pairs have coefficients of class $C^{q+2}$ on $N$;
the spaces $\mathcal P_l$ are understood in their Sobolev completions.
\begin{enumerate}[label=\textup{(\roman*)}]

\item Define the linear normal correction and reduced Hessian by
\begin{align}
 w_1(b)a&=-\mathfrak B_{b,Z}^{-1}\Pi_Z\mathfrak B_b\iota(a),
 &T_ba&=\iota(a)+w_1(b)a,
 \label{eq:linear-normal-correction}\\
 \mathfrak h_b[a,c]&=\langle\mathfrak B_bT_ba,T_bc\rangle,
 &a,c&\in K_+.
 \label{eq:reduced-hessian}
\end{align}
Then $w_1(b)a\in Z$ is linear in $a$ and analytic in $b$, and
\[
 \Pi_K(T_ba)_{\mathrm t}=a,\qquad
 \Pi_Z\mathfrak B_bT_ba=0,\qquad
 w_1(0)=0,\quad T_0a=\iota(a).
\]
If $a\in K_l$, then $T_ba\in\mathcal P_l$.
The form $\mathfrak h_b$ is symmetric and is identified with a
self-adjoint operator on $K_+$ using the fixed inner product. And we have
\[
 \|w_1(b)a\|_{H^q}\le C_q|b||a|.
\]

\item The test tensor
$r_b:=\Pi_{K_+}(j_b)_{\mathrm t}$ belongs to $K_2$ and satisfies
\begin{equation}
 \begin{gathered}
 r_0=v_*,\qquad |r_b-v_*|\le C|b|,\qquad
 \|T_br_b-j_b\|_{H^{q+2}}\le C_q|\mathrm De(b)|,\\
 |\mathfrak h_br_b|\le C|\mathrm De(b)|.
 \end{gathered}
 \label{eq:radial-schur}
\end{equation}

\item Define the quadratic correction by
\begin{equation}
 \begin{split}
 w_2(b,a)&=-\frac{1}{2}\mathfrak B_{b,Z}^{-1}\Pi_Z
       \mathrm D^2\mathcal G(H_b)[T_ba,T_ba],\\
 H^{[2]}(b,a)&=H_b+T_ba+w_2(b,a).
 \end{split}
 \label{eq:quadratic-normal-correction}
\end{equation}
Then $w_2(b,a)\in Z\cap\mathcal P_{\le L_1}$ is homogeneous
quadratic in $a$ and analytic in $b$ with coefficients of the prescribed
$C^{q+2}(N)$ regularity. And we have
\[
 \|w_2(b,a)\|_{H^q}\le C_q|a|^2.
\]
In any fixed Hermite bases, the coefficients on $N$ of $w_1(b)a$
and $w_2(b,a)$ have $C^q$ norms bounded by $C_q|b||a|$ and
$C_q|a|^2$, respectively. Moreover,
\begin{equation}
 H^{[2]}\in\mathcal S,\qquad
 \Pi_K(H^{[2]})_{\mathrm t}=b+a,\qquad
 \mathrm T_a^2[\Pi_Z\mathcal G(H^{[2]})]=0.
 \label{eq:normal-taylor-cancellation}
\end{equation}

\item Define the reduced energy by
\begin{equation}
 F_2(b,a)=e(b)+\frac{1}{2}\mathfrak h_b[a,a]
 +\frac{1}{6}\mathrm D^3\mathscr W_{H_b}
                   [T_ba,T_ba,T_ba].
 \label{eq:second-order-reduced-energy}
\end{equation}
It is analytic in $b$, polynomial of degree at most three in $a$,
and satisfies
\[
 F_2(b,a)=\mathrm T_a^3[\mathscr W(H^{[2]}(b,a))].
\]
In particular,
\[
 F_2(b,0)=e(b),\qquad \mathrm D_aF_2(b,0)=0,\qquad
 \mathrm D_a^2F_2(b,0)=\mathfrak h_b.
\]
Identifying derivatives with vectors in $K=K_N\oplus K_+$ by the
fixed inner product, one has the full gradient comparison
\begin{equation}
 \mathrm DF_2(b,a)
 =\mathrm T_a^2[\Pi_K\mathcal G(H^{[2]}(b,a))_{\mathrm t}]+O(|a|^3),
 \label{eq:finite-gradient-comparison}
\end{equation}
uniformly for $b$ in a smaller fixed neighborhood.
\end{enumerate}
\end{proposition}
\begin{proof}
Prescribe $q$ and apply Lemma~\ref{lem:moving-normal-inverse}
for the spaces $Z_l$ with $l\le L_1$ at a Sobolev order exceeding
$q+2+m/2$. Shrink the parameter neighborhood accordingly.
Products of coefficients on $N$ and Sobolev embedding then give
all the finite-order bounds below.

\textup{(i)} For $a\in K_+$, the linearized normal equation is
\[
 \mathfrak B_{b,Z}w_1(b)a=-\Pi_Z\mathfrak B_b\iota(a),
\]
which gives \eqref{eq:linear-normal-correction}. Since
$\mathfrak B_0\iota(a)=0$, the coefficient perturbation estimate
and the normal inverse give $\|w_1(b)a\|_{H^q}\le C_q|b||a|$.
Preservation of $\mathcal P_l$ gives $T_bK_l\subset\mathcal P_l$.
The remaining linear identities follow from $w_1(b)a\in Z$ and
the definition of $Z$. Symmetry of $\mathfrak B_b$ gives symmetry of
$\mathfrak h_b$. The inverse and the coefficients are analytic in
$b$ in these Sobolev spaces, giving the claimed finite-order
coefficient regularity.

\textup{(ii)} By Proposition~\ref{prop:moving-radial-hessian} and the projections
onto $\mathcal P_l$, $r_b\in K_2$, $r_0=v_*$, and $|r_b-v_*|\le C|b|$.
Also $j_b-\iota(r_b)\in Z$. Thus
$\delta_b:=T_br_b-j_b\in Z\cap\mathcal P_2$ satisfies
\[
 \mathfrak B_{b,Z}\delta_b=-\Pi_Z\mathfrak B_bj_b.
\]
Estimate~\eqref{eq:moving-radial-full} and the normal inverse give the third
bound in \eqref{eq:radial-schur}. For $c\in K_+$,
\[
 \langle\mathfrak B_b\delta_b,T_bc\rangle
 =\langle\delta_b,\mathfrak B_bT_bc\rangle=0,
\]
because $\Pi_Z\mathfrak B_bT_bc=0$. Consequently,
\[
 \mathfrak h_b[r_b,c]
 =\langle\mathfrak B_bj_b,T_bc\rangle,
\]
which proves $|\mathfrak h_br_b|\le C|\mathrm De(b)|$.

\textup{(iii)} The quadratic source in \eqref{eq:quadratic-normal-correction}
contains two factors from $T_ba\in\mathcal P_{\le J}$,
at most two spatial derivatives, and at most two explicit factors $x_i$.
This follows by differentiating
\eqref{eq:physical-constrained-gradient}: inverse metrics and density
denominators have coefficients depending only on $N$ at $H_b$, while the only
unbounded background coefficient is
$\mathrm df_E=(1/2)\sum_i x_i\,\mathrm dx_i$.
Its ordinary polynomial degree is therefore at most $2J+2$.
For a Euclidean--Euclidean component, the shifted index is two
larger than its Hermite degree; hence the entire source lies in
$\mathcal P_{\le 2J+4}=\mathcal P_{\le L_1}$. Lemma~\ref{lem:hermite-block-calculus}
and the normal inverse preserve this finite sum. The $\mathcal P_0$
component of the projected source lies in $Z_0$.
The inverse estimates give the stated Sobolev bounds for $w_2$.
Sobolev embedding on $N$ at higher fixed orders gives the coefficient
$C^q$ bounds for both corrections. Their homogeneity and analyticity
follow directly from their formulas.

The formal Taylor expansion is
\[
 \mathcal G(H^{[2]})=\mathcal G(H_b)
       +\mathfrak B_bT_ba+\mathfrak B_bw_2
 +\frac{1}{2}\mathrm D^2\mathcal G(H_b)
                  [T_ba,T_ba]+O_a(|a|^3).
\]
Its degree-zero normal term vanishes because
$\Pi_Z\mathcal G(H_b)=0$. The degree-one term vanishes by the
equation for $w_1$, and the degree-two terms cancel by the equation
for $w_2$. Both corrections belong to $Z$, so the slice condition
and kernel projection identity in \eqref{eq:normal-taylor-cancellation}
also follow.

\textup{(iv)} Since $\mathcal G(H_b)=\iota(\nabla e(b))$, its pairing with
$T_ba$ and with $w_2$ is zero. Moreover,
\[
 \mathrm D^2\mathscr W_{H_b}[T_ba,w_2]
 =\langle\mathfrak B_bT_ba,w_2\rangle=0.
\]
Thus the linear term vanishes, and $w_2$ contributes neither a
quadratic nor a cubic energy term. Hence the formal Taylor
polynomial $\mathrm T_a^3[\mathscr W(H^{[2]})]$ equals the expression
\eqref{eq:second-order-reduced-energy}.

Finally, for $\beta\in K_N$ and $\alpha\in K_+$,
\[
 \mathrm DH^{[2]}[\beta,\alpha]
 =\iota(\beta+\alpha)+Z_{\beta,\alpha},\qquad
 Z_{\beta,\alpha}\in Z.
\]
Indeed, derivatives of $v_0$, $w_1$, and $w_2$ take values in the
fixed normal space. The normal term pairs with
$\Pi_Z\mathcal G(H^{[2]})=O_a(|a|^3)$, so the formal chain rule gives
\[
 \mathrm D[\mathscr W(H^{[2]})][\beta,\alpha]
 =\langle\Pi_K\mathcal G(H^{[2]})_{\mathrm t},\beta+\alpha\rangle
       +O_a(|a|^3)(|\beta|+|\alpha|).
\]
Differentiating $\mathrm T_a^3$ in $a$ retains exactly the terms through
degree two in the $K_+$ component. Differentiating in $b$ retains
the derivative $\mathrm De(b)$ and the terms through degree three; its
cubic term is $O(|a|^3)$ on a smaller fixed $b$-ball.
This proves \eqref{eq:finite-gradient-comparison} in both components.
\end{proof}

\subsection{Cutoff construction and remainder estimates}
\label{sec:modulation-nonlinear}

The polynomial pair $H^{[2]}(b,a)$ need not define a positive metric
and density on the whole product, even for small parameters. We therefore
construct a cutoff family $\hat H_{b,a}$ with positive metric and density
satisfying the scalar constraints exactly, with errors flat in $a$ in
the weighted norms.

\begin{lemma}
\label{lem:positive-representatives}
Prescribe an integer $j\ge0$ and a finite number of parameter
derivatives.
On a sufficiently small parameter neighborhood depending on these
orders, there are coordinate
pairs $\hat H_{b,a}$, with $\hat H_{b,0}=H_b$, such that
\[
 \int_{\bar M}(\hat H_{b,a})_{\mathrm s}\;\mathrm d\bar\mu=0,
 \qquad
 \int_{\bar M}\vec x\,(\hat H_{b,a})_{\mathrm s}\;\mathrm d\bar\mu=0,
\]
\[
 \bar g+(\hat H_{b,a})_{\mathrm t}>0,
 \qquad 1+(\hat H_{b,a})_{\mathrm s}>0.
\]
The metric is uniformly equivalent to $\bar g$, and
$\hat X_{b,a}=\mathcal C(\hat H_{b,a})$ is complete, normalized,
and centered. They satisfy
\[
 \|\hat H_{b,a}-H_b\|_{C^{j+2}}\le C_j|a|^{1/2}.
\]
Moreover, for every $A>0$,
\begin{equation}
 \|\hat H_{b,a}-H^{[2]}(b,a)\|_{H^{j+2}}
 +\|\delta_{\bar X}(\hat H_{b,a})_{\mathrm t}\|_{H^{j+1}}
 +|\Pi_K(\hat H_{b,a})_{\mathrm t}-(b+a)|
 \le C_{j,A}|a|^A.
 \label{eq:positive-representative-flat}
\end{equation}
Here flat means $O(|a|^A)$ for every $A>0$ in the stated norms,
uniformly for $b$ in a smaller fixed neighborhood. These errors
remain flat after the prescribed $b$- and $a$-derivatives, and the
differentiated estimates extend at $a=0$.
The families agree on overlaps of their parameter neighborhoods.
\end{lemma}
\begin{proof}
Let $L_1$ be as in Proposition~\ref{prop:second-order-normal-reduction}.
Fix
\[
 0<\eta<\frac{1}{2(L_1+1)}
\]
and a smooth cutoff $0\le\chi\le1$ on $\mathbb R^d$, equal to one
on $\{|\vec x|\le1\}$ and zero on $\{|\vec x|\ge2\}$.
These choices are independent of the prescribed regularity orders.
Choose fixed smooth compactly supported functions
$\zeta_0,\ldots,\zeta_d$ on $\mathbb R^d$ such that
\[
 \begin{aligned}
 \int_{\bar M}\zeta_0\;\mathrm d\bar\mu&=1,&
 \int_{\bar M} x_i\zeta_0\;\mathrm d\bar\mu&=0,\\
 \int_{\bar M}\zeta_j\;\mathrm d\bar\mu&=0,&
 \int_{\bar M} x_i\zeta_j\;\mathrm d\bar\mu&=\delta_{ij}\quad(j\ge1).
 \end{aligned}
\]
They are obtained from a radial even bump and suitably normalized
multiples by $x_j$.
Apply Proposition~\ref{prop:second-order-normal-reduction} at a
finite order sufficient for the prescribed derivatives and the
$C^{j+2}$ estimate, shrinking the parameter neighborhood as needed.
Set $U_{b,a}=H^{[2]}(b,a)-H_b$. Then
$U_{b,a}\in\mathcal S\cap\mathcal P_{\le L_1}$, with coefficient
norms on $N$ bounded by $O(|a|)$ at all the orders required below. Its ordinary
polynomial degree is at most $L_1$. Since spatial differentiation
does not increase that degree, at each required order $k$,
\[
 |\bar\nabla^kU_{b,a}|\le C_k|a|\langle\vec x\rangle^{L_1}.
\]
For $0<|a|<1$, put $R_a=|a|^{-\eta}\ge1$ and
$\chi_a(\vec x)=\chi(\vec x/R_a)$. Define
\[
 c_0=\int_{\bar M}\chi_a(U_{b,a})_{\mathrm s}\;\mathrm d\bar\mu,
 \qquad c_i=\int_{\bar M} x_i\chi_a(U_{b,a})_{\mathrm s}\;\mathrm d\bar\mu
 \quad(1\le i\le d),
\]
\[
 \hat H_{b,a}
 =H_b+\chi_aU_{b,a}
       -\left(0,\sum_{\alpha=0}^dc_\alpha\zeta_\alpha\right).
\]
At $a=0$ set $\hat H_{b,0}=H_b$. Both $H_b$ and $U_{b,a}$
have zero scalar mean and first moments. Thus
\[
 c_0=-\int_{\bar M}(1-\chi_a)(U_{b,a})_{\mathrm s}\;\mathrm d\bar\mu,
 \qquad
 c_i=-\int_{\bar M} x_i(1-\chi_a)(U_{b,a})_{\mathrm s}\;\mathrm d\bar\mu.
\]
The prescribed moments of the $\zeta_\alpha$ give the two exact
scalar constraints for $\hat H_{b,a}$.

Every spatial derivative of $1-\chi_a$ is supported where
$|\vec x|\ge R_a$, and derivatives of positive order contribute
nonpositive powers of $R_a$. Gaussian integration of the polynomial
bounds therefore gives, for each required spatial order and fixed polynomial
weight, bounds of the form
\[
 C|a|R_a^M e^{-cR_a^2}
\]
for the weighted Sobolev norms of $(1-\chi_a)U_{b,a}$ and for
$|c_\alpha|$. Here $M$ is fixed by the orders in question and $c>0$.
Since $R_a=|a|^{-\eta}$, these bounds are $O(|a|^A)$ for every
$A>0$. The prescribed parameter derivatives introduce only fixed
powers of $|a|^{-1}$ and $R_a$, which are absorbed by the same
Gaussian exponential. Thus the tails and the moment corrections are
flat, with these derivatives, in every required polynomially weighted
Sobolev norm; these corrections have coefficients $c_\alpha$ and fixed
spatial supports.

We have
\begin{align*}
 \hat H_{b,a}-H^{[2]}(b,a)
   &=(\chi_a-1)U_{b,a}
       -\left(0,\sum_{\alpha=0}^dc_\alpha\zeta_\alpha\right),\\
 \delta_{\bar X}(\hat H_{b,a})_{\mathrm t}
   &=(U_{b,a})_{\mathrm t}(\bar\nabla\chi_a,\cdot),\\
 \Pi_K(\hat H_{b,a})_{\mathrm t}-(b+a)
   &=\Pi_K\bigl[(\chi_a-1)(U_{b,a})_{\mathrm t}\bigr].
\end{align*}
The second follows from
$\delta_{\bar X}(H_b)_{\mathrm t}
 =\delta_{\bar X}(U_{b,a})_{\mathrm t}=0$;
the scalar correction changes no tensor component. The last uses the
kernel projection identity in Proposition~\ref{prop:second-order-normal-reduction}.
The tail bounds and boundedness of $\Pi_K$ prove
\eqref{eq:positive-representative-flat}.

Leibniz's rule also gives
\[
 \|\chi_aU_{b,a}\|_{C^{j+2}}
 \le C_j|a|R_a^{L_1}
 =C_j|a|^{1-\eta L_1}\le C_j|a|^{1/2}.
\]
The flat scalar corrections satisfy the same estimate. At order zero,
smallness of $b$ and $a$ gives uniform metric equivalence to $\bar g$
and a uniform positive lower bound for the density. Completeness
follows from metric equivalence, and the exact scalar constraints
give normalization and centering. Equation~\eqref{eq:density-physical-comparison}
gives the corresponding bounded metric--potential estimates.

Extending the flat errors and their prescribed derivatives by zero at
$a=0$ gives the asserted extension. The normal inverses agree on
overlapping parameter neighborhoods, and the cutoff and correction
functions are fixed; hence the cutoff families agree there as well.
\end{proof}

\begin{lemma}
\label{lem:positive-cubic-residual}
For every prescribed integer $j\ge0$, choose the cutoff family
of Lemma~\ref{lem:positive-representatives} at sufficiently high
finite regularity. On a sufficiently small parameter neighborhood
depending on $j$,
\begin{equation}
 \|\Pi_Z\mathcal G(\hat H_{b,a})\|_{H^j}
 +|\Pi_K\mathcal G(\hat H_{b,a})_{\mathrm t}-\mathrm DF_2(b,a)|
 \le C_j|a|^3.
 \label{eq:positive-reduced-gradient}
\end{equation}
Here $\mathrm DF_2$ includes both parameter components and is
identified with a vector in $K_N\oplus K_+$ by the fixed inner product.
\end{lemma}
\begin{proof}
We prove the stronger estimate
\[
 \bigl\|\mathcal G(\hat H_{b,a})
       -\mathrm T_a^2[\mathcal G(H^{[2]}(b,a))]\bigr\|_{H^j}
 \le C_j|a|^3.
\]
Set
$V=\hat H_{b,a}-H_b$ and $H_t=H_b+tV$, $0\le t\le1$.
The metric and density remain positive along $H_t$, since both are
affine in $t$. Taylor's formula for the ambient extension in
\eqref{eq:physical-constrained-gradient} gives
\[
 \mathcal G(H_b+V)
 =\mathcal G(H_b)+\mathrm D\mathcal G(H_b)[V]
       +\frac{1}{2}\mathrm D^2\mathcal G(H_b)[V,V]
 +\frac{1}{2}\int_0^1(1-t)^2
           \mathrm D^3\mathcal G(H_t)[V,V,V]\;\mathrm dt.
\]
The differential here acts on $V$, which need
not belong to the metric slice.

The preceding construction gives, at every required fixed order,
\[
 |\bar\nabla^kV|\le C_k|a|\langle\vec x\rangle^{L_1},
 \qquad
 \|V\|_{C^{j+2}}\le C_j|a|^{1/2}.
\]
The scalar corrections obey the first bound because their coefficients
are flat. In the local expressions of
\eqref{eq:physical-constrained-gradient}, inverse metrics and density
denominators stay uniformly controlled along $H_t$, as do their
spatial derivatives at each fixed order. Every trilinear term in the
third differential contains three factors involving $V$ or its spatial
derivatives, and only fixed polynomial growth in $\vec x$ from the
background. Leibniz's rule thus bounds all derivatives through order
$j$ of these unprojected expressions by
$C_j|a|^3\langle\vec x\rangle^{M_j}$ for a fixed $M_j$.
Gaussian integration, followed by the fixed Sobolev projection bounds
in Lemma~\ref{lem:weighted-linear}, yields
\[
 \|\mathrm D^3\mathcal G(H_t)[V,V,V]\|_{H^j}
       \le C_j|a|^3
 \qquad(0\le t\le1).
\]
The bounded $O(|a|^{1/2})$ norm controls the coefficients and
positivity; the polynomial $O(|a|)$ bound on each of the three
factors involving $V$ supplies the full cubic power.

Write
\[
 V=T_ba+w_2(b,a)+E_{b,a},
 \qquad E_{b,a}=\hat H_{b,a}-H^{[2]}(b,a).
\]
The polynomially weighted tail estimates in
Lemma~\ref{lem:positive-representatives} make all
terms containing $E_{b,a}$ flat in $H^j$. The norms of the coefficients
on $N$ of $T_ba$ and $w_2$ are respectively $O(|a|)$ and
$O(|a|^2)$. The mixed quadratic term is therefore $O_{H^j}(|a|^3)$,
and the term with two copies of $w_2$ is $O_{H^j}(|a|^4)$.
Hence the terms through quadratic order are exactly
$\mathrm T_a^2[\mathcal G(H^{[2]}(b,a))]$, proving the stronger
estimate above. Equations~\eqref{eq:normal-taylor-cancellation}
and~\eqref{eq:finite-gradient-comparison} give
\eqref{eq:positive-reduced-gradient}.
\end{proof}

\section{A rigidity inequality for cylinders with compact Einstein factors}
\label{sec:Einstein-rigidity}
Throughout this section, assume $2\Ric_{g_N}=g_N$.
Then $f_N$ is constant and $v_*=u_*g_N\in K_G\subset K_2$.
\subsection{The quadratic obstruction on \texorpdfstring{$K_+$}{K+}}
\label{sec:Gaussian-quadratic}

All pairings use the fixed product metric and probability measure.

\begin{theorem}
\label{thm:radial-obstruction}
For $h=\sum_{j\ge0}h_j\in K$, we have
\begin{equation}
 \int_{\bar M}\left\langle
 \mathrm D^2(\Phi\circ\mathcal C)_0[\iota(h),\iota(h)],v_*
 \right\rangle\;\mathrm d\bar\mu
 =-\sum_{j\ge1}j^2\|h_j\|_{L^2}^2,
 \qquad \|v_*\|_{L^2}^2=8md.
 \label{eq:full-radial-obstruction}
\end{equation}
\end{theorem}
\begin{proof}
Along the fixed-density path
$X_{t,s}=(g_{t,s},f_{t,s})
=\mathcal C(t\iota(h)+s\iota(v_*))$,
\eqref{eq:physical-entropy-variation} gives
\[
 \partial_s\mathcal W(X_{t,s})
 =\int_{\bar M}\langle\Phi(X_{t,s}),v_*\rangle_{g_{t,s}}
       \;\mathrm d\bar\mu.
\]
Differentiating twice in $t$, the derivatives of the metric pairing
multiply $\Phi(\bar X)$ or
$\mathrm D(\Phi\circ\mathcal C)_0\iota(h)$, both of which vanish.
Hence
\[
 \mathrm D^3\mathscr W_0[\iota(h),\iota(h),\iota(v_*)]
 =\int_{\bar M}\left\langle
 \mathrm D^2(\Phi\circ\mathcal C)_0[\iota(h),\iota(h)],v_*
 \right\rangle\;\mathrm d\bar\mu.
\]
We compute this cubic using explicit metric--potential paths.
For a metric $G$ on $N$, set
\[
 f_G^0:=f_N+\frac{1}{2}\log\frac{\det G}{\det g_N},\qquad
 E_N(G):=\int_N
 (\mathrm R_G+|\nabla f_G^0|_G^2+f_G^0-m)\;\mathrm d\mu_N.
\]
For $G=G(p,\vec{x})$, with $p\in N$, use the pair $X=(g,f)$ given by
\[
 g=g_E+G,\qquad f=f_E+f_G^0.
\]
Its measure is $\mathrm d\bar\mu$. Put $B_i:=G^{-1}\partial_iG$.
We use indices $a,b$ tangent to $N$ and Euclidean indices $i,j$, with
$\partial_i=\partial_{x_i}$. Besides the Christoffel symbols of $G$ on $N$, the additional
connection coefficients are
\[
 \Gamma^a_{ib}=\frac{(B_i)^a_b}{2},\qquad
 \Gamma^i_{ab}=-\frac{(\partial_iG)_{ab}}{2}.
\]
Substitution into the Ricci formula gives
\begin{align*}
 (\Ric_g)_{ij}
 &=-\frac{1}{2}\partial_j\Tr B_i
       -\frac{1}{4}\Tr(B_iB_j),\\
 (\Ric_g)_{ab}
 &=(\Ric_G)_{ab}-\frac{1}{2}\sum_i\partial_i^2G_{ab}
       -\frac{1}{4}\sum_i(\Tr B_i)\partial_iG_{ab}
       +\frac{1}{2}\sum_i
       (\partial_iG\,G^{-1}\partial_iG)_{ab}.
\end{align*}
Since
\[
 \Tr(G^{-1}\partial_i^2G)
 =\partial_i\Tr B_i+\Tr(B_i^2),
\]
the trace over $TN$ is
\[
 G^{ab}(\Ric_g)_{ab}
 =\mathrm R_G-\frac{1}{2}\sum_i\partial_i\Tr B_i
       -\frac{1}{4}\sum_i(\Tr B_i)^2.
\]
Adding the Euclidean trace yields
\[
 \mathrm R_g=\mathrm R_G-\sum_i\partial_i\Tr B_i
 -\frac{1}{4}\sum_i\Tr(B_i^2)
 -\frac{1}{4}\sum_i(\Tr B_i)^2.
\]
Since $\partial_if=\frac{x_i}{2}+\frac{\Tr B_i}{2}$, the trace-square terms cancel
in $\mathrm R_g+|\nabla f|^2$. Gaussian integration by parts cancels
the remaining derivative and mixed terms. Thus
\begin{equation}
 \mathcal W(X)=\int_{\mathbb R^d}E_N(G(\vec{x}))\;\mathrm d\gamma_d
 -\frac{1}{4}\int_{\bar M}\sum_i
 \Tr(G^{-1}\partial_iG\,G^{-1}\partial_iG)\;\mathrm d\bar\mu.
 \label{eq:entropy-block-splitting}
\end{equation}

Write $E_N=A+B-m$, where $B(G)=\int_N f_G^0\;\mathrm d\mu_N$.
Constant scaling gives
\[
 E_N(cG)=c^{-1}A(G)+B(G)+\frac{m}{2}\log c-m.
\]
At $g_N$, a divergence-free tensor $\eta$ on $N$ satisfies
\[
 \mathrm D^2E_N[\eta,\eta]=\frac{1}{2}\int_N
 \langle\mathcal L_N\eta,\eta\rangle\;\mathrm d\mu_N,
 \qquad \mathrm D^2B[\eta,\eta]=-\frac{1}{2}\|\eta\|_{L^2}^2.
\]
Differentiate twice in $t$ and once in $c$ at $(c,t)=(1,0)$ in
$E_N(cg_N+t\eta)$. This yields
\begin{equation}
 \mathrm D^3E_N[g_N,\eta,\eta]
 =-\frac{3}{2}\int_N\langle\mathcal L_N\eta,\eta\rangle\;\mathrm d\mu_N
 -\frac{1}{2}\|\eta\|_{L^2}^2.
 \label{eq:compact-scaling-third}
\end{equation}
Apply \eqref{eq:entropy-block-splitting} to
$G=g_N+th+su_*g_N$. Expanding its two inverse matrices, the
$t,t,s$ derivative of the last integral is
\[
 \int_{\bar M}\left(u_*|\nabla_{\vec{x}}h|^2+
 2\sum_i\partial_iu_*\langle h,\partial_ih\rangle\right)
 \;\mathrm d\bar\mu.
\]
Consequently the third entropy derivative is
\[
 \int_{\bar M}\left(-\frac{3}{2}u_*\langle\mathcal L_Nh,h\rangle
 -\frac{1}{2}u_*|h|^2+u_*|\nabla_{\vec{x}}h|^2
 +2\sum_i\partial_iu_*\langle h,\partial_ih\rangle\right)
 \;\mathrm d\bar\mu.
\]
Use $\mathcal L_Nh=-\Delta_{f_E}h$ and $\Delta_{f_E}u_*=-u_*$.
Two integrations by parts reduce this expression to
\begin{equation}
 -\frac{1}{4}\int_{\bar M}u_*|h|^2\;\mathrm d\bar\mu
 -\frac{1}{2}\int_{\bar M}u_*|\nabla_{\vec{x}}h|^2\;\mathrm d\bar\mu.
 \label{eq:radial-Hermite-pairing}
\end{equation}
By Lemma~\ref{lem:hermite-identities}, on $\mathcal H_j$ the degree-preserving
part of multiplication by $|\vec{x}|^2$ is $(2d+4j)\Id$.
Euclidean derivatives lower Hermite degree by one,
and $\int_{\bar M}|\nabla_{\vec{x}}h_j|^2\;\mathrm d\bar\mu=(j/2)\|h_j\|_{L^2}^2$.
It follows that
\[
 \int_{\bar M}u_*|h_j|^2\;\mathrm d\bar\mu=4j\|h_j\|_{L^2}^2,
 \qquad
 \int_{\bar M}u_*|\nabla_{\vec{x}}h_j|^2\;\mathrm d\bar\mu
 =2j(j-1)\|h_j\|_{L^2}^2.
\]
For different $j$, the tensor coefficients on $N$ belong to orthogonal
eigenspaces of the self-adjoint operator $\mathcal L_N$. This
orthogonality eliminates all cross terms in
\eqref{eq:radial-Hermite-pairing}, including the Hermite degrees
separated by two that multiplication by $|\vec{x}|^2$ can couple.
The third derivative is therefore $-\sum_jj^2\|h_j\|_{L^2}^2$,
which proves the first identity in \eqref{eq:full-radial-obstruction}.

For comparison with the second variation in $(g,f)$, the chain rule gives
\[
 \mathrm D^2(\Phi\circ\mathcal C)_0
       [\iota(h),\iota(h)]
 =\mathrm D^2\Phi_{\bar X}
       [\mathrm D\mathcal C_0\iota(h),
        \mathrm D\mathcal C_0\iota(h)]
       +\frac{1}{2}\bar\nabla^2|h|^2.
\]
The last term pairs to zero with every element of $K$, since
$\delta_{\bar X}K=0$. Thus both radial pairings, and their
$K$-projections, agree.

To justify the polynomial variations, first multiply them by a smooth
cutoff equal to one on $|\vec{x}|\le R$ and supported on
$|\vec{x}|\le2R$. The metrics and densities along these paths are then positive for small parameters.
The errors in the differentiated identities and integrations by parts
are bounded by a fixed polynomial times a Gaussian tail and vanish
as $R\to\infty$. Finally,
$\|v_*\|_{L^2}^2=m\int_{\mathbb R^d} u_*^2\;\mathrm d\gamma_d=8md$.
\end{proof}

\begin{corollary}
\label{cor:full-quadratic-coercivity}
For
\[
 Q(h):=\frac{1}{2}\Pi_K
 \mathrm D^2(\Phi\circ\mathcal C)_0[\iota(h),\iota(h)],
 \qquad h=b+a\in K_N\oplus K_+,
\]
we have
\begin{equation}
 \|Q(h)\|_{L^2}
 \ge\frac{\sum_{j\ge1}j^2\|a_j\|_{L^2}^2}{2\sqrt{8md}}
 \ge\frac{\|a\|_{L^2}^2}{2\sqrt{8md}}.
 \label{eq:full-quadratic-coercivity}
\end{equation}
\end{corollary}
\begin{proof}
Since $v_*\in K$, the projection does not change its pairing.
Apply Cauchy--Schwarz to \eqref{eq:full-radial-obstruction}.
\end{proof}

\subsection{Proof of the rigidity inequality}
\label{sec:local-mixed-rigidity}

\begin{lemma}
\label{lem:finite-closure}
Suppose an analytic function $F(b,a)$ satisfies $F(b,0)=e(b)$,
$\mathrm D_aF(b,0)=0$, and a bounded continuous vector $r_b$ satisfies
\[
 |\mathrm D_a^2F(b,0)r_b|\le C|\mathrm De(b)|,
 \qquad \mathrm D_a^3F(0,0)[r_0,a,a]\le-c|a|^2.
\]
Then, on a smaller neighborhood,
$|a|^2+|\mathrm De(b)|\le C|\mathrm DF(b,a)|$.
\end{lemma}
\begin{proof}
Since $\mathrm D_aF(b,0)=0$, differentiating in $b$ gives
$\mathrm D_a\mathrm D_bF(b,0)=0$. Taylor expansion of
$\mathrm D_bF(b,a)$ and $\mathrm D_aF(b,a)[r_b]$ therefore gives,
on a smaller neighborhood,
\begin{align*}
 |\mathrm De(b)|&\le|\mathrm DF(b,a)|+C|a|^2,\\
 \frac{c}{4}|a|^2
 &\le C|\mathrm DF(b,a)|+C|\mathrm De(b)||a|.
\end{align*}
For the second estimate, continuity preserves the strict negativity of
$\mathrm D_a^3F(b,0)[r_b,a,a]$ for $a\ne0$, while the Hessian term is bounded by
$C|\mathrm De(b)||a|$.
Substituting the first estimate gives
\[
 |a|^2\le C\bigl((1+|a|)|\mathrm DF(b,a)|+|a|^3\bigr).
\]
Choose $|a|$ small enough to absorb the cubic term. This bounds
$|a|^2$ by $C|\mathrm DF(b,a)|$, and the first estimate then gives
the same bound for $|\mathrm De(b)|$.
\end{proof}

Apply Lemma~\ref{lem:finite-closure} to the reduced energy in
Proposition~\ref{prop:second-order-normal-reduction}. That proposition gives
\[
 F_2(b,0)=e(b),\qquad \mathrm D_aF_2(b,0)=0,\qquad
 |\mathrm D_a^2F_2(b,0)r_b|
       =|\mathfrak h_br_b|\le C|\mathrm De(b)|.
\]
The field $r_b$ is bounded and continuous. Since
$T_0=\iota$ and $r_0=v_*$, the cubic energy formula in
Proposition~\ref{prop:second-order-normal-reduction} and the
fixed-density identity in the proof of
Theorem~\ref{thm:radial-obstruction} give
\begin{equation}
 \mathrm D_a^3F_2(0,0)[r_0,a,a]
       =-\sum_{j\ge1}j^2\|a_j\|_{L^2}^2\le-|a|^2.
 \label{eq:finite-radial-third}
\end{equation}
Lemma~\ref{lem:finite-closure} therefore gives
\begin{equation}
 |a|^2+|\mathrm De(b)|\le C|\mathrm DF_2(b,a)|.
 \label{eq:finite-closed}
\end{equation}

For an actual pair $X=\mathcal C(H)$ with $H=(h,\sigma)$ satisfying
\[
 \int_{\bar M}\sigma\;\mathrm d\bar\mu=0,\qquad
 \int_{\bar M}\vec x\,\sigma\;\mathrm d\bar\mu=0,
\]
define, for every integer $j\ge0$,
\begin{equation}
 \D_j(X)=\|\Phi(X)\|_{H^j}
 +\|\delta_{\bar X}h\|_{H^{j+1}},
 \qquad \D=\D_0.
 \label{eq:augmented-defect-full}
\end{equation}
For pairs satisfying $\delta_{\bar X}h=0$, $\D=\|\Phi\|_{L^2}$.
Set $(b,a)=\Pi_Kh$ and $\Xi=H-\hat H_{b,a}$.

\begin{lemma}
\label{lem:remainder-projection}
Fix an integer $j\ge0$ and choose the cutoff family $\hat H_{b,a}$ at the
corresponding finite regularity. If
$H\in H^{j+2}$ satisfies the scalar constraints above and $(b,a)$ lies
in the parameter neighborhood of Lemma~\ref{lem:positive-representatives},
then, for every $A>0$,
\begin{equation}
 \|(I-\Pi_Z)\Xi\|_{H^{j+2}}
 \le C_j\|\delta_{\bar X}h\|_{H^{j+1}}
       +C_{j,A}|a|^A.
 \label{eq:constraint-projection-bound}
\end{equation}
\end{lemma}
\begin{proof}
Both coordinate endpoints satisfy the scalar constraints, so the scalar
component of $(I-\Pi_{\mathcal S})\Xi$ vanishes. The Hodge projection
formula is therefore
\[
 (I-\Pi_Z)\Xi=
 \left(\delta_{\bar g}^*\bar{\mathcal P}^{-1}
       \delta_{\bar X}\Xi_{\mathrm t}+\Pi_K\Xi_{\mathrm t},\,0\right).
\]
Here $\delta_{\bar X}\Xi_{\mathrm t}$ is $L^2(\mathrm d\bar\mu)$-orthogonal
to the one-forms dual to Killing fields by adjointness. Lemma~\ref{lem:weighted-linear} bounds the inverse, and
Lemma~\ref{lem:positive-representatives} gives
\[
 \delta_{\bar X}\Xi_{\mathrm t}
       =\delta_{\bar X}h+O_{H^{j+1}}(|a|^A),\qquad
 \Pi_K\Xi_{\mathrm t}=O(|a|^A).
\]
Finite-dimensional norm equivalence on $K$ now proves
\eqref{eq:constraint-projection-bound}.
\end{proof}

\begin{theorem}
\label{thm:local-rigidity-input}
Assume the finite-order obstruction on $N$, and set
$\rho=\max\{2,\ell\}$. There is a fixed bounded $C^2$
neighborhood of zero in metric--density coordinates such that every
$H=(h,\sigma)$ in this neighborhood with
\[
 \bar g+h>0,\qquad 1+\sigma>0,\qquad
 \int_{\bar M}\sigma\;\mathrm d\bar\mu=0,\qquad
 \int_{\bar M}\vec x\,\sigma\;\mathrm d\bar\mu=0
\]
defines a complete normalized centered pair $X=\mathcal C(H)$ satisfying
\begin{equation}
 |a|^2+|b|^\ell+\|\Xi\|_{H^2}\le C\D(X),\qquad
 \|H\|_{H^2}\le C\D(X)^{1/\rho}.
 \label{eq:full-local-rigidity}
\end{equation}
Here $(b,a)=\Pi_Kh$, $\D$ is \eqref{eq:augmented-defect-full}, and
$\Xi=H-\hat H_{b,a}$. For every integer $j\ge0$, in a sufficiently
small bounded $C^{j+2}$ neighborhood contained in the preceding
$C^2$ neighborhood,
\begin{equation}
 \|\Xi\|_{H^{j+2}}\le C_j\D_j(X).
 \label{eq:full-high-normal}
\end{equation}
This neighborhood and its constant may depend on $j$.
When $K_N=0$, the $b$ terms are absent and $\rho=2$.
\end{theorem}
\begin{proof}
Finite-dimensional norm equivalence and the probability normalization give
\[
 |b|+|a|\le C\|h\|_{L^2}\le C\|H\|_{C^0}.
\]
Thus smallness of $H$ places $(b,a)$ in the fixed parameter neighborhoods
required for the order-zero argument. Lemma~\ref{lem:positive-representatives}
and the estimate for $H_b$ give
\[
 \|\hat H_{b,a}\|_{C^2}\le C|b|+C|a|^{1/2}.
\]
After shrinking the original neighborhood, the affine segment from
$\hat H_{b,a}$ to $H$ is therefore sufficiently small in $C^2$.
The associated metric and density are positive along this segment,
and both scalar constraints hold exactly.

Integrate \eqref{eq:gradient-coefficient-comparison} along the segment
to obtain
\[
 \mathcal G(H)-\mathcal G(\hat H_{b,a})
 =\mathrm D\mathcal G(0)\Xi+R,
 \qquad \|R\|_{L^2}\le C\epsilon\|\Xi\|_{H^2}.
\]
On $Z$, the derivative is $\mathfrak B_0$, with inverse
\eqref{eq:full-normal-inverse}. For every ambient coordinate pair $U$,
Lemma~\ref{lem:constrained-entropy-gradient} and self-adjointness give
\[
 \Pi_K\bigl(\mathrm D\mathcal G(0)U\bigr)_{\mathrm t}
 =\frac12\Pi_K\mathcal L_{\bar f}U_{\mathrm t}=0.
\]
Project the gradient comparison onto $Z$ and $\iota(K)$.
Lemma~\ref{lem:remainder-projection}, with $j=0$ and $A=3$, controls
$(I-\Pi_Z)\Xi$. Together with the cubic residual and gradient comparison
\eqref{eq:positive-reduced-gradient} and the defect bound
\eqref{eq:gradient-by-defect}, this yields
\[
 \|\Xi\|_{H^2}+|\mathrm DF_2(b,a)|
 \le C\bigl(\D(X)+|a|^3\bigr)+C\epsilon\|\Xi\|_{H^2}.
\]
Absorb the last term by choosing $\epsilon$ small enough.
Equation~\eqref{eq:finite-closed} then yields
$|a|^2+|\mathrm De(b)|\le C(\D+|a|^3)$.
Absorb the cubic term on a fixed small neighborhood.
Assumption~\ref{ass:compact-rigidity} gives
$|b|^\ell\le C|\mathrm De(b)|$, proving the first inequality.
The cutoff family satisfies
$\|\hat H_{b,a}\|_{H^2}\le C(|b|+|a|)$, so for $\D\le1$
the second follows with $\rho=\max\{2,\ell\}$.

For each fixed integer $j\ge0$, choose the cutoff family $\hat H_{b,a}$ at the
required finite regularity and shrink the bounded $C^{j+2}$ neighborhood.
The affine segment is then small at this order. Repeat the normal
comparison, using Lemma~\ref{lem:remainder-projection} with $A=3$, to obtain
\[
 \|\Xi\|_{H^{j+2}}
 \le C_j(\D_j+|a|^3)+C_j\epsilon_j\|\Xi\|_{H^{j+2}}.
\]
Absorb the last term by choosing $\epsilon_j$ sufficiently small. The
order-zero estimate already proved gives
$|a|^3\le C\D^{3/2}\le C\D_j$ on a smaller neighborhood.
This proves
\eqref{eq:full-high-normal}.
\end{proof}

\section{Global rigidity}
\label{sec:global-rigidity}
We follow the contraction--extension construction of
\cite[Section~5]{LiZhang2023}. In its contraction proof, the local
rigidity inequality is applied at (5.23) and the normal remainder is
estimated at (5.25). Here Theorem~\ref{thm:local-rigidity-input} supplies
these two inputs for the
full kernel and arbitrary finite obstruction order on $N$.

\subsection{Contraction and extension}
\label{sec:contraction-extension}
Throughout this subsection, $N^m$ is a connected compact Einstein manifold,
$m\ge2$, with $2\Ric_{g_N}=g_N$, satisfying
Assumption~\ref{ass:compact-rigidity}, and $d\ge1$.
For a probability-normalized shrinker put
$F=f-\boldsymbol\mu(g,1)$. Thus
$\mathrm R+|\nabla F|^2=F$, and
\begin{equation}
 r=2\sqrt{F},\qquad
 \bar r=2\sqrt{\bar f-\Theta_N}=\sqrt{|\vec x|^2+2m}.
 \label{eq:static-radial-functions}
\end{equation}
The potential $F$ uses the normalization of \cite{LiZhang2023}.

Put $\rho=\max\{2,\ell\}$. The geometric extension proof below fixes
$\kappa_{\mathrm{ext}}>0$ and an upper bound $\theta_*\le1/10$
independently of the comparison rates. Choose once
\begin{equation}\label{eq:comparison-rates}
 \alpha=\min\left\{\frac{1}{32\rho},\frac{\kappa_{\mathrm{ext}}}{4}\right\},
 \qquad \beta=\frac{\alpha}{8},\qquad
 \beta'=\frac{9\beta}{10},\qquad 0<\theta\le\theta_*.
\end{equation}
These choices satisfy
\begin{equation}\label{eq:alpha-contraction}
 \rho\alpha<\beta'+\frac18,\qquad
 \alpha<\min\left\{\frac18,\kappa_{\mathrm{ext}}\right\},
 \qquad \beta(1+\theta)^2<\frac\alpha2.
\end{equation}
Fix one sufficiently large finite derivative order $k_*$, as specified
in the contraction proof, and constants $D_j\ge1$ and nonnegative
integer powers $A_j$ dominating the extension bounds for
$1\le j\le k_*$. The order, constants, powers, and smallness thresholds
are fixed before iteration. The obstruction order on $N$ enters
through $\rho$; the additional modes in $K_+$ have fixed finite Hermite
degree and contribute only fixed polynomial factors.

\begin{definition}[Strong and weak comparison conditions]
\label{def:star-dagger}
Both conditions use a diffeomorphism $\chi:\Omega\to\{r<L\}$,
where $\Omega\subset\bar M$ is connected and relatively compact, with
smooth boundary in $\{L-10^{-2}<\bar r<L+10^{-2}\}$, and
\[
 \{\bar r<L-10^{-2}\}\subset\Omega
 \subset\{\bar r<L+10^{-2}\}.
\]
For a tensor or function set
\begin{equation}\label{eq:static-bracket-norm}
 [T]_j=\sum_{i=0}^j|\bar\nabla^iT|_{\bar g}.
\end{equation}
Condition $(\star_L)$ means
\[
 [\chi^*g-\bar g]_5+[\chi^*f-\bar f]_5
 \le e^{\alpha(\bar r^2-L^2)}
\]
on $\Omega$, with the left side also at most a fixed sufficiently
small $\epsilon_*>0$ throughout $\Omega$.
Condition $(\dagger_L)$ means
\[
 \begin{split}
 |\chi^*g-\bar g|+|\chi^*f-\bar f|&\le e^{-\beta L^2},\\
 [\chi^*g-\bar g]_j+[\chi^*f-\bar f]_j&\le D_jL^{A_j}
 \quad(1\le j\le k_*).
 \end{split}
\]
\end{definition}

\begin{theorem}[Contraction]\label{thm:higher-order-contraction}
There are fixed $c_0,L_1>0$ such that
$(\dagger_L)\Longrightarrow(\star_{L-c_0})$ for $L\ge L_1$.
\end{theorem}
\begin{proof}
Fix an integer $q$ with $q\beta'>1/8$, and choose $k_*$ sufficiently
large to perform these gauge corrections and obtain the estimates
through order $n+8$ used below.
The cutoff and interpolation of
\cite[proof of Theorem~5.2, (5.16)--(5.18)]{LiZhang2023} give, through
the orders needed for $q$ gauge steps, an ordinary derivative bound
$CL^P\epsilon_L$, where $\epsilon_L=e^{-\beta'L^2}$.
Here and below $P$ denotes a fixed power, enlarged when necessary.
The central identity $\boldsymbol\mu=f-\mathrm R-|\nabla f|^2$
gives a uniform entropy lower bound. The volume and
potential estimates on sublevel sets of $r$ \cite[Lemmas~3.1--3.2]{LiZhang2023} then give
uniform Gaussian tails.

Localize the pair across a collar of fixed width and set it equal to
the reference pair outside the chart. Its total mass differs from one
by $O(L^Pe^{-(L-c)^2/4})$. Adding the logarithm of this mass to its
potential makes the mass exactly one without changing $\Phi$ or the
reference metric divergence. Denote this normalized pair by $X$ and
write $\mathrm d\mu_X=(1+\sigma)\,\mathrm d\bar\mu$. The cutoff and
normalization retain $\|\sigma\|_{C^0}\le CL^P\epsilon_L$.
Indeed, the density ratio is the product of $e^{-(f-\bar f)}$ and the
relative volume density; both are uniformly close to one by the
preceding estimates.

Choose a radial cutoff $\eta$ supported in the comparison region and
equal to one on a smaller interior. For the vector parameter
$t\in\R^d$, let $\psi_t$ be the time-one flow of
$V_t=\eta T_t$, where $T_t=\sum_i t_i\partial_{x_i}$.
Change of variables and differentiation of the inverse flow give
\begin{align*}
 m(t)&:=\int_{\bar M}\vec x\;\mathrm d\mu_{\psi_t^*X}
      =\int_{\bar M}(\vec x\circ\psi_t^{-1})\;\mathrm d\mu_X,\\
 \mathrm Dm(0)&=-\left(\int_{\bar M}\eta\;\mathrm d\mu_X\right)I.
\end{align*}
The central mass lower bound and the uniform cutoff derivative bounds
make the quantitative inverse function theorem uniform in $L$.
Since $m(0)=\int_{\bar M}\vec x\,\sigma\;\mathrm d\bar\mu$, it gives a unique
small $t_*$ with
\[
 m(t_*)=0,\qquad |t_*|\le C|m(0)|
 \le C\|\sigma\|_{C^0}\le CL^P\epsilon_L.
\]
Pullback preserves total mass, so the recentered scalar density
coordinate has exactly zero mean and first moments.

On nested interior regions $U'\Subset U$ whose flow trajectories stay
in $\{\eta=1\}$, one has $\psi_t(p,\vec x)=(p,\vec x+t)$ and
$\psi_t^*\bar g=\bar g$. For the actual metric $g=\bar g+h$,
\[
 \begin{gathered}
 \psi_t^*g-g=h(p,\vec x+t)-h(p,\vec x),\\
\|\psi_t^*g-g\|_{C^k(U')}
\le C|t|\,\|h\|_{C^{k+1}(U)}.
 \end{gathered}
\]
Thus its initial interior change is $O(L^P\epsilon_L^2)$, after
enlarging $P$. The linear background contribution
\[
L_{\eta T_t}\bar g
=\mathrm d\eta\otimes T_t^\flat+T_t^\flat\otimes\mathrm d\eta
\]
is confined to the outer cutoff collar and gains a Gaussian tail in
weighted estimates. The potential change is bounded by $CL|t|$.
After each metric gauge step the same construction centers its new
moment exactly; the interior metric error gains the additional small
factor from $h$, preserving the gauge improvement below.

Alternate the localized metric gauge of
\cite[Theorem~3.13 and Lemma~5.1, especially (5.15)]{LiZhang2023}
with this exact recentering. The gauge inverse is taken $L^2(\mathrm d\bar\mu)$-orthogonal to
Killing fields. Cutting off its source requires a projection onto the
finite-dimensional space of Killing fields; integration by parts places its coefficients in the outer
collar, of size $CL^P\epsilon_Le^{-(L-c_i)^2/4}$. They can be removed
with fixed compactly supported dual fields. If $s_i$ bounds the
interior metric divergence after step $i$, the localized flow estimates
give
\[
 s_{i+1}\le CL^P\epsilon_Ls_i
       +CL^P\epsilon_Le^{-(L-c_i)^2/4},
 \qquad s_0\le CL^P\epsilon_L.
\]
Indeed, the moment change of the intermediate metric flow is at most
$CL^P(\epsilon_Ls_i+(s_i+\epsilon_L)e^{-(L-c_i)^2/4})$:
the uncut gauge field is $L^2(\mathrm d\bar\mu)$-orthogonal to the
translation fields, its
removed tail has this bound, and replacing the reference density
costs $CL^P\epsilon_L$. Exact recentering therefore has the same size;
its interior metric-divergence error gains another factor
$CL^P\epsilon_L$. This proves the displayed recurrence also for the
alternating construction.

After these fixed $q$ steps, let $X_L=\mathcal C(H_L)$ be the last
centered pair. The flows preserve mass and no further cutoff is made.
Consequently $H_L=(h,\sigma)$ has exactly zero scalar mean and first
moments. The soliton defect vanishes on the retained interior; the
divergence on the interior is $CL^P\epsilon_L^{q+1}$ plus the tail above.
The collar estimates and $q\beta'>1/8$ give
\begin{equation}\label{eq:localized-defect}
 \D_{n+8}(X_L)\le CL^P
 \exp\left[-\beta'L^2-\frac{(L-c_g)^2}{8}\right]=:\delta_L,
\end{equation}
where $c_g$ is a fixed total collar width. Scalar moment corrections
use weighted $L^1$ tails with exponent $1/4$, whereas this weighted
$L^2$ estimate uses exponent $1/8$.

Applying Theorem~\ref{thm:local-rigidity-input} to $X_L$, with
$\Xi=H_L-\hat H_{b,a}$, gives
\begin{equation}\label{eq:loc-input}
 |a|^2+|b|^\ell\le C\delta_L,
 \qquad \|\Xi\|_{H^{n+10}}\le C\delta_L.
\end{equation}
Lemma~\ref{lem:positive-representatives} gives
the polynomial comparison bound
\[
[\hat H_{b,a}]_5\le C(1+|\vec x|)^P(|a|+|b|)
            \le C(1+|\vec x|)^P\delta_L^{1/\rho}.
\]
On reference balls of radius $c(1+|\vec x|)^{-1}$, with $c$ below
the product injectivity radius, the Gaussian density varies by a
bounded factor relative to $e^{-|\vec x|^2/4}$. Since
$n+10>5+n/2$, rescaled Sobolev embedding gives
\[
[\Xi]_5\le C(1+|\vec x|)^P e^{|\vec x|^2/8}
                  \|\Xi\|_{H^{n+10}}
         \le C(1+|\vec x|)^P e^{|\vec x|^2/8}\delta_L.
\]
Adding these bounds and applying \eqref{eq:density-physical-comparison}
yields
\begin{equation}\label{eq:pointwise-separated}
 [h]_5+[f-\bar f]_5
 \le CL^P\left(\delta_L^{1/\rho}
             +e^{\bar r^2/8}\delta_L\right)
 \quad\text{on }\{\bar r\le L-c_g-2\}.
\end{equation}
Since $|a|$ is exponentially small, the cutoff radius used to define
$\hat H_{b,a}$, a negative power of $|a|$, lies outside this comparison chart.

Choose $c_0>c_g+3$. The inequality
$\rho\alpha<\beta'+1/8$ absorbs the polynomial in the kernel term.
For the normal term, division by
$e^{\alpha(\bar r^2-(L-c_0)^2)}$ gives at most
\[
 CL^P\exp\left[-\beta'L^2+
 \frac{(L-c_0)^2-(L-c_g)^2}{8}\right]
\]
on the retained interior, since $\alpha<1/8$. This tends to zero.
Restore the initial mass-normalizing constant in the potential; its
$CL^Pe^{-(L-c)^2/4}$ error is smaller than the same comparison weight.
As at the end of \cite[proof of Theorem~5.2]{LiZhang2023}, restrict the
composed chart to the entire target sublevel $\{r<L-c_0\}$.
The remaining collar and the radial derivative of the geometric potential
give a boundary graph at distance $o(L^{-1})$ from the model level,
with the required connected domain and core. The error tends uniformly
to zero there, which also gives the $\epsilon_*$ condition.
\end{proof}

\begin{theorem}[Extension]\label{thm:higher-order-extension}
There are fixed $\theta,L_2>0$ such that
$(\star_L)\Longrightarrow(\dagger_{(1+\theta)L})$ for $L\ge L_2$.
\end{theorem}
\begin{proof}
The curvature extension and radial matching are
\cite[proof of Theorem~5.3, Steps~1--2, (5.30)--(5.36)]{LiZhang2023}.
Their hypotheses are bounded geometry of the model and small $C^5$
comparison on the interior collar. Pseudolocality extends curvature
control to a fixed larger radius, and the potential flow supplies the
radial identification. Glue the resulting chart to the model outside
a slightly smaller extension, obtaining a complete weighted manifold
$(\bar M,g',F')$ which agrees with the given shrinker on
$\{r'<(1+\tau/2)L\}$ and with the reference cylinder on
$\{r'>(1+\tau)L\}$, where $r'=2\sqrt{F'}$ and $\tau>0$ is fixed
and sufficiently small. Its geometry is uniformly close to the model
in the sense of \cite[Remark~A.2]{LiZhang2023}; radial gluing also gives
uniform bounds for $\Ric_{g'}$ and $\nabla^2F'$.

The scalar spectrum of $-\Delta_{\bar f}$ has precisely $d$ eigenvalues
at $1/2$. Indeed,
$\lambda_1(-\Delta_N)\ge m/(2(m-1))>1/2$, and the next positive
Gaussian level is $1$. This gap persists under the small global gluing
perturbation. Apply \cite[proof of Theorem~5.3, Step~3, (5.37)--(5.39)]{LiZhang2023} to its
$d$-dimensional spectral cluster. The cutoff coordinate test functions
give an error $CL^Pe^{-\alpha L^2}$; Bochner's formula includes the
gluing-annulus error as in the Bochner identity between (5.38) and (5.39)
of that proof. Choosing $\tau$ small first
makes that error at most $CL^Pe^{-L^2/8}$, and gives
\begin{equation}\label{eq:eigen-small}
 |\lambda_i-\tfrac12|+\|\nabla^2v_i\|_{L^2(e^{-F'}\mathrm dV_{g'})}^2
 \le CL^Pe^{-\alpha L^2}.
\end{equation}
We follow the growth step in the proof of Theorem~5.3, Step~4, of
\cite{LiZhang2023}. Proposition~A.1, (A.4)--(A.5), of that paper gives
$U_{v_i}=D_{v_i}/I_{v_i}\le4\lambda_i$ and, with $\theta=0$,
polynomial growth of $I_{v_i}$ on the glued manifold. On the interior
region where the glued pair agrees with the given shrinker,
$R_1=R_2=\mathrm R\ge0$, so (A.12) gives
$I_v'\le2D_v/r$ above a fixed radius. Interior elliptic estimates
therefore give polynomial bounds for $v_i$ and their derivatives on
smaller sublevel sets within this region. Take
\[
 R=(1+3\tau/8)L,\qquad S=(1+\tau/4)L,
 \qquad \kappa_{\mathrm{ext}}
 =\frac{(1+3\tau/8)^2-(1+\tau/4)^2}{6}.
\]
For $T=\nabla^2v_i$, the exact shrinker commutator gives
$\langle\Delta_{F'}T,T\rangle\ge-\Lambda|T|^2$ up to $R$, with
fixed $\Lambda>0$. Lemma~A.3, (A.22), of \cite{LiZhang2023},
or the frequency identities in
\cite[Lemmas~4.7, 4.15 and proof of Proposition~4.16]{ColdingMinicozzi2025}
give, with $Q=4\Lambda+2\eta$, $0<\eta<\min\{1,\Lambda\}$, and fixed large $r_1$,
\begin{equation}\label{eq:localized-hessian-growth}
 I_T(t)\le C(t/r_1)^Q
 \left(I_T(r_1)+D_T(R)_+e^{-(R^2-t^2)/6}\right),
 \qquad r_1\le t\le S.
\end{equation}
Here, at a regular level $r'=r$, the frequency quantities are
\begin{align*}
 I_T(r)&=r^{1-n}\int_{\{r'=r\}}|T|^2|\nabla r'|\;\mathrm dA_{g'},\\
 D_T(r)&=\frac{r^{2-n}}{2}\int_{\{r'=r\}}
 \left\langle\nabla|T|^2,\frac{\nabla r'}{|\nabla r'|}\right\rangle
 \;\mathrm dA_{g'}.
\end{align*}
To see the boundary term, set $B=D_T-(2\Lambda+\eta)I_T$; the cited differential inequality is
$B'\ge rB/3$ where $B>0$. Integrate it backwards from $R$, then
integrate $I_T'\le2D_T/r$ with the factor $r^{-Q}$.
The preceding scalar bound gives $D_T(R)_+\le CL^P$, while
\eqref{eq:eigen-small} gives $I_T(r_1)\le CL^Pe^{-\alpha L^2}$
at a fixed central level. Thus \eqref{eq:localized-hessian-growth}
has rate $\min\{\alpha,\kappa_{\mathrm{ext}}\}$ for the squared
Hessian norm. No shrinker commutator is used in the gluing annulus.

Since $\alpha<\kappa_{\mathrm{ext}}$, the localized estimate
\eqref{eq:localized-hessian-growth} and the normalization of the scalar functions in
\cite[proof of Theorem~5.3, Step~4, preceding (5.43)]{LiZhang2023} yield, on a slightly smaller
sublevel and through the fixed required orders,
\begin{equation}\label{eq:almost-splitting-new}
 |\langle\nabla u_i,\nabla u_j\rangle-\delta_{ij}|
 +[\nabla^2u_i]_k
 +|2\langle\nabla F',\nabla u_i\rangle-u_i|
 \le CL^Pe^{-\alpha L^2/2}=:E_L.
\end{equation}

Let $N_0$ be the component of $\{u=0\}$ near the original copy of $N$.
Identify it with
$N$ by the original central comparison
chart, as in \cite[proof of Theorem~5.3, Step~5, preceding (5.44)]{LiZhang2023}.
The restrictions of $u_i$ to the original central copy have small
tangential Hessians; the elliptic estimate on $N$ modulo constants
makes their tangential gradients small. Their normal gradients form
an almost orthonormal frame, and the last term in
\eqref{eq:almost-splitting-new} makes their central values small.
Thus $N_0$ is a small normal graph with metric and geometric potential
within $CL^PE_L$ of $g_N$ and $m/2$. The slice on $N$ and elliptic
regularity in \cite[(5.46)--(5.47)]{LiZhang2023} give coordinate
regularity.

The almost-parallel flow construction of
\cite[proof of Theorem~5.3, Steps~5--6, (5.44)--(5.47)]{LiZhang2023} gives the enlarged chart.
One may use the dual fields
$Y_i=\sum_j(Q^{-1})_{ij}\nabla u_j$,
$Q_{ij}=\langle\nabla u_i,\nabla u_j\rangle$, and their simultaneous
flow with parameter $\vec x$, so $u=t\vec x$ along the flow.
Reversing this flow gives injectivity. Its potential estimate places
the source boundary outside the desired target sublevel, so restriction
maps onto that entire sublevel. Indeed every sublevel component has an
interior minimum with $F=\mathrm R\le n/2$, hence
meets the already covered connected core. The radial derivative gives
the boundary band and core containment in Definition~\ref{def:star-dagger}.
For a fixed $\theta_*\le\min\{\tau/8,1/10\}$, these constructions give
\begin{equation}\label{eq:extension-endpoint}
 |\chi^*g-\bar g|+|\chi^*f-\bar f|
 \le CL^Pe^{-\alpha L^2/2},
\end{equation}
on the chart to radius $(1+\theta)L$, $0<\theta\le\theta_*$.
Here the entropy shift is restored using its
$O(e^{-\alpha L^2})$ bound from the original central chart.
Curvature bounds, the slice on $N$, and differentiation of the flow
give $D_jL^{A_j}$ bounds through the fixed order $k_*$. All their
constants and powers are independent of $L$. The last inequality in
\eqref{eq:alpha-contraction} absorbs the fixed polynomial factor and
proves the weak comparison.
\end{proof}

\begin{proof}[Proof of Theorem~\ref{thm:intro-rigidity}]
Proposition~\ref{prop:compact-isolation} gives
Assumption~\ref{ass:compact-rigidity} for some finite order.
The initial smooth comparison in \cite[proof of Theorem~5.4]{LiZhang2023}
shows that sufficient pointed closeness gives $(\dagger_{L_0})$ at
any prescribed sufficiently large radius $L_0$.
Theorems~\ref{thm:higher-order-contraction} and
\ref{thm:higher-order-extension} then give comparisons at radii
\[
 L_{j+1}=(1+\theta)(L_j-c_0).
\]
These radii tend to infinity when $L_0$ exceeds the fixed iteration
threshold. The exhaustion argument in the same reference yields a
global isometry preserving the normalized potentials.
\end{proof}

\subsection{Finite quotients}
\label{sec:finite-quotients}
Let $\Gamma$ be a finite group acting freely and isometrically on
$\bar M=N\times\R^d$.
After translating the common fixed point of its Euclidean action to
the origin, fix the marking
$\gamma(p,\vec x)=(I_\gamma p,A_\gamma\vec x)$, $A_\gamma\in O(d)$.
For $\pi:\bar M\to\bar M_\Gamma=\bar M/\Gamma$, the normalization is
\begin{equation}
 \pi^*\bar g_\Gamma=\bar g,\qquad
 \pi^*\bar f_\Gamma=\bar f-\log|\Gamma|,\qquad
 \Theta_\Gamma=\Theta_N-\log|\Gamma|.
 \label{eq:quotient-normalization}
\end{equation}
Thus the geometric potential and radius pull back to
$\bar f-\Theta_N$ and $\bar r$. A marked comparison chart is a
diffeomorphism $\hat\chi_L:\Omega_\Gamma\to\{r<L\}$ with the
core and boundary conditions of Definition~\ref{def:star-dagger}.
Set $\Omega=\pi^{-1}(\Omega_\Gamma)$.
The lift $\chi_L=\hat\chi_L\circ\pi|_\Omega$ is a principal $\Gamma$-cover
of the target sublevel. We impose the same strong or weak bounds on
\[
 \tilde g=\chi_L^*g,\qquad
 \tilde f=\chi_L^*f+\log|\Gamma|,
\]
relative to $(\bar g,\bar f)$, and denote them by
$(\star_L^\Gamma)$ and $(\dagger_L^\Gamma)$.
Subtracting $\boldsymbol\mu(g,1)+\log|\Gamma|$ from
$\tilde f$ recovers $\chi_L^*F$, so the
target geometric radius also agrees after lifting.
The lifted density is $|\Gamma|^{-1}$ times the pulled-back target
density, so invariant integrals agree with those on the target sublevel.
All comparison constants are fixed before iteration and may depend on
the marked action.

\begin{theorem}
\label{thm:quotient-rigidity}
Let $(N^m,g_N)$ be a smooth closed connected Einstein manifold with
$\Ric_{g_N}=\frac{g_N}{2}$, $m>0$, rigid among normalized
compact Ricci shrinkers.
For the fixed free finite action above and every
$\bar p\in\operatorname*{argmin}\bar f_\Gamma$, there exists
$\epsilon>0$ such that any connected complete normalized
$n$-dimensional Ricci shrinker $(M,g,f,p)$ with
$p\in\operatorname*{argmin}f$ and
\[
 d_{\mathrm{pGH}}((M,g,p),(\bar M_\Gamma,\bar g_\Gamma,\bar p))
 <\epsilon
\]
is isometric to $(\bar M_\Gamma,\bar g_\Gamma)$.
\end{theorem}
\begin{proof}
Proposition~\ref{prop:compact-isolation} gives
Assumption~\ref{ass:compact-rigidity} for some finite order.
We use the quotient construction of
\cite[Theorems~5.6--5.7]{LiZhang2023}, retaining the marking throughout.
The fixed projections and reconstruction commute with the action.
The family $H_b$, normal corrections, and gauge inverse are
equivariant by uniqueness on their fixed complements. Choose radial
cutoffs and average the finite-dimensional right inverses used for
corrections enforcing the scalar constraints and Killing corrections. The invariant radial direction
in Corollary~\ref{cor:full-quadratic-coercivity} gives its estimate on
the full $K_+^\Gamma$.

For an invariant normalized density, the first-moment vector belongs
to $(\R^d)^\Gamma$. Perform the exact recentering in this subspace;
the derivative of the first-moment map remains a uniformly invertible scalar
multiple of the identity. If $(\R^d)^\Gamma=0$, the moment already
vanishes. The alternating gauge and recentering construction is
therefore equivariant and gives
$(\dagger_L^\Gamma)\Rightarrow(\star_{L-c_0}^\Gamma)$ with the
same scalar constraints and defect bound as in
Theorem~\ref{thm:higher-order-contraction}.

For extension, begin with this finite cover of the comparison sublevel.
The local curvature extension and the potential flow extend it across
the enlarged radial collar, as in
\cite[proof of Theorem~5.7, (5.49) and the following paragraph]{LiZhang2023}.
On that collar $|\nabla F|^2=F-\mathrm R>0$;
the flow therefore extends the deck action and the cover of the entire
sublevel. The extension takes place on the cover of the comparison sublevel.
Apply the glued-model construction of
Theorem~\ref{thm:higher-order-extension} on the extended cover, with
invariant gluing data, as in \cite[proof of Theorem~5.7, paragraph following (5.49)]{LiZhang2023}.

The scalar spectral projection onto the entire $d$-dimensional cluster
near $1/2$ commutes with $\Gamma$. Project the cutoff coordinate vector
and normalize its Gram matrix equivariantly. The resulting vector $u$
satisfies $u(\gamma q)=A_\gamma u(q)$. The component $N_0$ of
$\{u=0\}$ near the original copy of $N$ is an equivariant small
normal graph over that copy. Its identification and
slice adjustment on $N$ can be chosen equivariantly.
The Gram matrix satisfies
$Q(\gamma q)=A_\gamma Q(q)A_\gamma^T$, so the simultaneous flow
generated by the dual fields $Y_i$ in the proof of
Theorem~\ref{thm:higher-order-extension} is equivariant. The resulting chart descends to the marked quotient,
retaining injectivity, surjectivity onto the entire sublevel set,
and the comparison bounds established there.
The extension estimate of Theorem~\ref{thm:higher-order-extension}
therefore gives
$(\star_L^\Gamma)\Rightarrow(\dagger_{(1+\theta)L}^\Gamma)$.

The initial smooth comparison argument in
\cite[proof of Theorem~5.4]{LiZhang2023} applies to the quotient. The limiting geometric potential lifts to
$m/2+|\vec x|^2/4$, and its total Gaussian mass gives the entropy
limit $\Theta_\Gamma$; hence the normalized potentials converge as
well. Apply the same radius recurrence and exhaustion to embeddings
of the marked quotient. The limit is injective by the local inverse function theorem
applied to the comparison embeddings, and proper by the potential
identity. It is consequently a global weighted isometry with this
specified quotient. In dimension three the result is the
round-cylinder quotient case of \cite[Theorem~1.2]{LiZhang2023}.
\end{proof}

\clearpage
\phantomsection
\addcontentsline{toc}{section}{References}
\begingroup
\fontsize{9.5}{11.2}\selectfont
\let\originalthebibliography\thebibliography
\renewcommand{\thebibliography}[1]{%
  \originalthebibliography{#1}%
  \setlength{\itemsep}{1.5pt}%
  \setlength{\parsep}{0pt}%
  \setlength{\parskip}{0pt}%
}
\bibliographystyle{rigidityamsalpha}
\bibliography{RigidityRicciShrinkersII_V29}
\endgroup

\vskip10pt

\noindent
Hanbing Fang, Department of Mathematics, The Chinese University of Hong
Kong, Shatin, N.T., Hong Kong; E-mail:
\href{mailto:hanbingfang@cuhk.edu.hk}{hanbingfang@cuhk.edu.hk}.

\smallskip

\noindent
Yu Li, Institute of Geometry and Physics, University of Science and Technology of China, No. 96 Jinzhai Road, Hefei, Anhui Province, 230026, China; Hefei National Laboratory, No. 5099 West Wangjiang Road, Hefei, Anhui Province, 230088, China; E-mail:
\href{mailto:yuli21@ustc.edu.cn}{yuli21@ustc.edu.cn}.

\end{document}